\documentclass[13pt,letterpaper,reqno]{article}
\usepackage{epsf}
\usepackage{amsmath}
\usepackage{amsfonts}
\usepackage{amssymb}
\usepackage{amstext}
\usepackage{amsbsy}
\usepackage{amsthm}
\usepackage{amscd}
\usepackage[all]{xy}

\newtheorem{theorem}{Theorem}[section]

\newtheorem{lem}[theorem]{Lemma}
\newtheorem{prop}[theorem]{Proposition}
\newtheorem{dfn}[theorem]{Definition}
\newtheorem{rem}[theorem]{Remark}
\newtheorem{cor}[theorem]{Corollary}
\newtheorem{cl}[theorem]{Claim}

\newcommand{\w}{\widetilde}

\newcommand{\bpf}{\noindent {\em Proof.  }}
\newcommand{\epf}{\qed \vspace{+10pt}}

\newcommand{\Zee}{\mathbb{Z}}

\newcommand{\Cee}{\mathbb{C}}

\newcommand{\eeight}{\left ( E_8 \times E_8 \right ) \rtimes \Zee_2}
\newcommand{\spint}{{\rm Spin}(32)/\Zee_2}

\newcommand{\hee}{{\rm H} \oplus {\rm E}_8 \oplus {\rm E}_8 }
\begin{document}

\title{Modular Invariants for Lattice Polarized K3 Surfaces}
\author{Adrian Clingher
\thanks{
Department of Mathematics, Stanford University, Stanford, CA 94305. {\bf e-mail:} {\it clingher@math.stanford.edu}
}
\and 
Charles F. Doran\thanks{
Department of Mathematics, University of Washington, Seattle, WA 98195. {\bf e-mail:} {\it doran@math.washington.edu}
}
} 
\maketitle
\begin{center} 
\abstract{\noindent We study the class of complex algebraic K3 surfaces admitting an embedding of 
$\hee$ inside the N\'{e}ron-Severi lattice. These special K3 surfaces are classified by a pair of 
modular invariants, in the same manner that elliptic curves over $\mathbb{C}$ are classified by 
the ${\rm J}$-invariant. Via the canonical Shioda-Inose structure we construct a geometric 
correspondence relating K3 surfaces of the above type with abelian surfaces realized as cartesian 
products of two elliptic curves. We then use this correspondence to determine explicit 
formulas for the modular invariants.}  
\end{center}
%
%
%
%
%
%
\section{Introduction}
Let ${\rm X}$ be an algebraic K3 surface over the field of complex numbers. The $\Zee$-module obtained as the image of the first Chern class map:
$$ c_1 \colon {\rm H}^1({\rm X}, \mathcal{O}^*_{{\rm X}}) \ \rightarrow \ {\rm H}^2({\rm X}, \Zee), $$
when endowed with the bilinear pairing induced by the intersection form on ${\rm H}^2({\rm X}, \Zee)$, 
forms an even lattice. By Lefschetz' Theorem on $(1,1)$ classes, this is precisely the 
N\'{e}ron-Severi lattice ${\rm NS}({\rm X})$ of the surface ${\rm X}$, namely the group of isomorphism classes 
of divisors modulo homological equivalence. Furthermore, according to the Hodge Index Theorem, ${\rm NS}({\rm X})$ is 
an indefinite lattice of rank $1 \leq {\rm p}_{{\rm X}} \leq 20$ and signature of type $(1, {\rm p}_{{\rm X}}-1)$. 
\par In \cite{dolgachev}, Dolgachev formulated the notion of a lattice polarization of a K3 surface. If
${\rm M}$ is an even lattice of signature $(1,r)$ with $r\geq 0$, then an ${\bf M}$-{\bf polarization} on 
${\rm X}$ is, by definition, a primitive lattice embedding:
\begin{equation}
\label{inoseemb}
i \colon {\rm M} \ \hookrightarrow \ {\rm NS}(X) 
\end{equation}  
such that the image $i({\rm M})$ contains a pseudo-ample class. A coarse moduli space $\mathcal{M}_{{\rm M}} $ 
can be defined for equivalence classes of pairs $({\rm X},i)$ of ${\rm M}$-polarized K3 surfaces and an appropriate version of the Global Torelli Theorem holds.  
\par The focus of this paper is on K3 surfaces which admit a polarization by the unique 
even unimodular lattice of signature $(1,17)$. This particular lattice can be realized effectively as 
the orthogonal direct sum  
$$ {\rm M} \ = \ {\rm H} \oplus {\rm E}_8 \oplus {\rm E}_8 $$
where ${\rm H}$ is the standard rank-two hyperbolic lattice and ${\rm E}_8$ is the unique even, negative-definite 
and unimodular lattice of rank eight. Note that not all algebraic K3 surfaces admit such an ${\rm M}$-polarization. In fact, the presence of such a structure imposes severe constraints on the geometry of ${\rm X}$. In particular, the Picard rank ${\rm p}_{{\rm X}}$ has to be $18, 19$ or $20$.   
\par A standard observation on the Hodge theory of this special class of K3 surfaces is that the polarized 
Hodge structure of an ${\rm M}$-polarized K3 surface $({\rm X},i)$ is identical with the polarized Hodge 
structure of an abelian surface ${\rm A}={\rm E}_1 \times {\rm E}_2$ realized as a cartesian product of 
two elliptic curves. Since both types of surfaces involved admit appropriate versions of the Torelli theorem, 
Hodge theory implies a well-defined correspondence:
\begin{equation}
\label{mirrormap}
({\rm X},i) \ \leftrightarrow \ {\rm E}_1 \times {\rm E}_2
\end{equation}
giving rise to a canonical analytic isomorphism between the corresponding moduli spaces on the 
two sides. By employing a modern point of view from the frontier of algebraic geometry with string 
theory, one can regard $(\ref{mirrormap})$ as a Hodge-theoretic {\bf duality map}, a correspondence that relates two 
seemingly different types of surfaces sharing similar Hodge-theoretic information\footnote{In fact, 
the use of this terminology for $(\ref{mirrormap})$ is quite natural. The 
identification of Hodge structures given by $(\ref{mirrormap})$ is a particular case of a more 
general Hodge-theoretic phenomenon which, surprisingly, was predicted by 
physics. In string theory this relationship is known as the F-Theory/Heterotic String Duality 
in eight dimensions. We refer the reader to Section \ref{finalcomments} for a brief discussion of 
this aspect.}. Our point in this work is that the resemblance of the two Hodge structures involved 
in the duality correspondence $(\ref{mirrormap})$ is not fortuitous, but rather is merely a consequence 
of a quite interesting geometric relationship. 
\begin{theorem}
\label{main1}
Let $({\rm X},i)$ be an ${\rm M}$-polarized K3 surface. 
\begin{itemize}
\item [(a)] The surface ${\rm X}$ possesses a canonical involution $\beta$ defining 
a Shioda-Inose structure. 
\item [(b)] The minimal resolution of ${\rm X}/\beta$ is a new K3 surface ${\rm Y}$ endowed 
with a canonical Kummer structure. This structure realizes ${\rm Y}$ as the Kummer surface 
${\rm Km}({\rm E}_1 \times {\rm E}_2)$ associated to an abelian surface ${\rm A}$ canonically 
represented as a cartesian product of two elliptic curves. The elliptic curves ${\rm E}_1$ and ${\rm E}_2$ 
are unique, up to permutation.
\item [(c)] The construction induces a canonical Hodge isomorphism between the 
${\rm M}$-polarized Hodge structure of ${\rm X}$ and the natural ${\rm H}$-polarized Hodge 
structure of the abelian surface ${\rm A}= {\rm E}_1 \times {\rm E}_2$. 
\end{itemize}
\end{theorem}
\noindent Section $\ref{partone}$ of the paper is devoted entirely to proving the above theorem. 
\par In the second part of the paper we describe an application of the geometric transform 
outlined above. One important feature of the special class of ${\rm M}$-polarized K3 surfaces 
is that such polarized pairs $({\rm X},i)$ turn out to be completely classified by two 
{\bf modular invariants} $\pi, \sigma \in \Cee $, much in the same as way elliptic curves over the field of complex numbers are classified by the J-invariant. However, the two modular invariants $\pi$ and $\sigma$ are not 
geometric in origin. They are defined Hodge-theoretically, and the result leading to the classification 
is a consequence of the appropriate version of the Global Torelli Theorem for lattice polarized 
K3 surfaces. However, in the context of the duality map $(\ref{mirrormap})$, the two invariants can be 
seen as the standard symmetric functions on the J-invariants of the dual elliptic curves:
\begin{equation}
\sigma \ = \ {\rm J}({\rm E}_1) + {\rm J}({\rm E}_2), \ \ \ \pi \ = \ {\rm J}({\rm E}_1) \cdot {\rm J}({\rm E}_2).
\end{equation}  
This interpretation suggests that the modular invariants of an ${\rm M}$-polarized 
K3 surface can be computed by determining the two elliptic curves that appear on the 
right-side of $(\ref{mirrormap})$.
\par Explicit ${\rm M}$-polarized K3 surfaces can be constructed by various geometrical 
procedures. One such method, introduced in 1977 by Inose \cite{inose}, constructs a 
two-parameter family ${\rm X}({\rm a}, {\rm b}) $ of ${\rm M}$-polarized K3 surfaces\footnote{
An equivalent two-parameter family is known in the physics literature as the Morrison-Vafa family 
\cite{morrisonvafa}.} by 
taking minimal resolutions of the projective quartics in 
$\mathbb{P}^3$ associated with the special equations:  
\begin{equation}
\label{inoseform111}
y^2zw - 4x^3z + 3{\rm a}xzw^2 - \frac{1}{2} \left ( z^2w^2 + w^4 \right ) + {\rm b} zw^3 \ = \ 0, \ \ \ 
{\rm a}, {\rm b} \in \Cee . \ 
\end{equation}
In fact, as we will see shortly, this family covers 
all possibilities. Every ${\rm M}$-polarized K3 surface can be realized as 
${\rm X}({\rm a}, {\rm b})$ for some ${\rm a}, {\rm b} \in \Cee $.  One can
regard the Inose quartic $(\ref{inoseform111})$ as a normal form of an 
${\rm M}$-polarized K3 surface.
\par In the second part of the paper, we use the geometric transform of Theorem $\ref{main1}$ to explicitly 
describe the J-invariants of the two elliptic curves ${\rm E}_1$ and ${\rm E}_2$ associated to the 
Inose surface ${\rm X}(a,b)$. 
\begin{theorem}
\label{main2}
The J-invariants ${\rm J}({\rm E}_1)$ and ${\rm J}({\rm E}_2)$ of the two elliptic curves associated 
to ${\rm X}(a,b)$ by the transform of Theorem $\ref{main1}$ are the two solutions of the 
quadratic equation:
$$ {\rm x}^2 - \left ( a^3-b^2+1 \right ) {\rm x} + a^3 \ = \ 0.$$ 
\end{theorem}
\noindent As pointed out earlier, as a consequence of the above theorem, one obtains explicit 
formulas for the two modular invariants of the Inose surface ${\rm X}(a,b)$. 
\begin{cor}
\label{mainth}
The modular invariants of the Inose surface ${\rm X}({\rm a},{\rm b})$ are given by:
\begin{equation}
\label{modulareq}
\pi = {\rm a}^3 , \ \ \ \  \sigma  = {\rm a}^3 - {\rm b}^2 + 1.
\end{equation} 
\end{cor}
\noindent The power of the geometric transformation underlying the duality map $(\ref{mirrormap})$ is fully revealed 
by the proof of Theorem $\ref{main2}$. In the absence of such a geometric argument, in order to prove that statement one would be forced to undertake long and very complex computations of the periods of the quartic $(\ref{inoseform111})$\footnote{Such an approach via period computations have been taken in the physics literature \cite{billo, curio}.}.  
\par It seems that the geometric transformation described by Theorem $\ref{main1}$ is a particular case 
of a more general phenomenon. Evidence for this is provided by an analysis of the slightly more general 
case of K3 surfaces polarized by the rank $17$ lattice ${\rm H} \oplus {\rm E}_8 \oplus {\rm E}_7$. 
Surfaces in this class still admit a canonical Shioda-Inose structure. This leads to a correspondence 
between these special K3 surfaces and jacobians of smooth genus-two curves. These results will be 
described in a forthcoming paper. 
\par A very interesting alternative arithmetic approach to the above questions, by Elkies and Kumar, 
has been communicated to the authors \cite{kumar}. 
\section{Hodge Structures for ${\rm M}$-polarized K3 Surfaces}  
\label{modularsection}
Let $({\rm X},i)$ be an ${\rm M}$-polarized K3 surface. Denote by $\omega \in {\rm H}^2({\rm X}, \Cee) $ the 
class of a non-zero holomorphic two-form on ${\rm X}$. This class is unique, up to multiplication by 
a non-zero scalar. The Hodge structure of ${\rm X}$ is then essentially given by the decomposition:
$$ {\rm H}^2({\rm X}, \Cee) \ = \ {\rm H}^{2,0}({\rm X}) \oplus {\rm H}^{1,1}({\rm X}) \oplus {\rm H}^{0,2}({\rm X}) $$
where ${\rm H}^{2,0}({\rm X}) = \Cee \cdot \omega $, 
${\rm H}^{0,2}({\rm X}) = \Cee \cdot \bar{\omega} $ 
and ${\rm H}^{1,1}({\rm X}) = \{ \omega, \bar{\omega} \} ^{\perp} $. Since the lattice $i({\rm M})$ is 
generated by classes associated to algebraic cycles, it follows that one has an embedding:
$$ i({\rm M}) \ \subset \  {\rm H}^{1,1}({\rm X}) \cap {\rm H}^2({\rm X}, \Zee). $$  
By standard lattice theory (see, for example, the exposition in \cite{nikulin}), the orthogonal complement ${\rm N}$ of 
$i({\rm M}) $ in ${\rm H}^2({\rm X},\Zee)$ is an even, unimodular sublattice of signature $(2,2)$. 
Hence the lattice ${\rm N}$ is isometric to the orthogonal direct sum ${\rm H} \oplus {\rm H}$ of 
two rank-two hyperbolic lattices. 
One can therefore choose a basis 
$$\mathcal{B} = \{ x_1,x_2,y_1,y_2 \}$$ 
of ${\rm N}$ with intersection matrix:
$$
\left ( 
\begin{array}{cccc}
0 & 0 & 1 & 0 \\
0 & 0 & 0 & 1 \\
1 & 0 & 0 & 0 \\
0 & 1 & 0 & 0 \\
\end{array}
\right ). 
$$  
It follows that $ \omega $ belongs to ${\rm N} \otimes \Cee $. Moreover, since the elements of the basis $ \mathcal{B}$ 
are isotropic, the class $\omega$ has non-zero intersection with any one of them. The class $\omega$ is therefore uniquely defined as soon as one imposes the normalization condition $(\omega, y_2)=1$.
\par Let us also note that the basis $ \mathcal{B} $ can be chosen such that the isomorphism of real vector spaces:
\begin{equation}
\label{orientrev}
(\omega, \cdot) \colon \langle x_1, x_2  \rangle  \otimes \mathbb{R} \ \rightarrow \ \Cee 
\end{equation}
is orientation reversing. Then, as discussed in \cite{clingher3}, the Hodge-Riemann bilinear relations imply that the normalized 
period class can be written:
\begin{equation}
\label{odescomp}
\omega \ = \tau x_1 + x_2 + u y_1 + (-\tau u) y_2  
\end{equation}
where $\tau, u$ are uniquely defined (but depending on the choice of basis $\mathcal{B}$) elements of the complex upper half-plane $\mathbb{H}$. 
\begin{dfn}
The {\bf modular invariants} of the ${\rm M}$-polarized K3 surface $({\rm X}, i)$ are, by definition:
\begin{equation}
\label{moddef}
\sigma({\rm X},i) := \ {\rm J}(\tau) + {\rm J}(u), 
\end{equation}
$$
\pi({\rm X},i) : = {\rm J}(\tau) \cdot {\rm J}(u)
$$
where ${\rm J}$ is the classical elliptic modular function\footnote{The function ${\rm J}$ is normalized such that 
${\rm J}(i)=1$ and ${\rm J}(e^{\frac{2 \pi i}{3}})=0$.}.
\end{dfn}
\noindent Let us make two observations justifying the importance of the numbers defined above. 
Firstly, the numbers $\sigma({\rm X},i)$ and $\pi({\rm X},i)$ do not depend on the choice of basis 
$\mathcal{B}$. That is because any new choice of basis $\mathcal{B}'$ can be related to 
$\mathcal{B}$ by an integral isometry $\varphi $ of the lattice ${\rm N}$ that preserves 
the spinor norm. That is $\mathcal{B}' = \varphi(\mathcal{B})$. But then, as shown for 
example in Section 2 of \cite{hosono2}, the group $O^+({\rm N})$ of such integral isometries is 
naturally isomorphic to the semi-direct product:
$$ \left ( {\rm PSL}(2,\Zee) \times {\rm PSL}(2,\Zee) \right ) \rtimes \Zee/  2 \Zee $$ 
with the generator of $\Zee/ 2 \Zee$ acting on  ${\rm PSL}(2,\Zee) \times {\rm PSL}(2,\Zee)$ by exchanging the 
two sides. This clearly proves that such a modification does not affect $\sigma({\rm X},i)$ and $\pi({\rm X},i)$. 
\par Secondly, thanks to a lattice polarized version of the Global Torelli Theorem for 
K3 surfaces \cite{dolgachev}, the numbers $\sigma({\rm X},i)$ and $\pi({\rm X},i)$ fully classify the 
polarized pairs $({\rm X},i)$ up to isomorphism. Simply put, this result says that there 
exists a two-dimensional complex analytic space $ \mathcal{M}_{{\rm M}} $ realizing a coarse moduli space 
for ${\rm M}$-polarized K3 surfaces and that the period 
map to the classifying space of polarized Hodge structures is an isomorphism of analytic spaces.
\begin{equation}
\label{hodgeident}
 \mathcal{M}_{{\rm M}} \ \stackrel{}{\longrightarrow} \ 
\left ( {\rm PSL}(2,\Zee) \times {\rm PSL}(2,\Zee) \right ) \rtimes \Zee / 2 \Zee \ \backslash \ 
\left ( \mathbb{H} \times \mathbb{H} \right ). 
\end{equation}
From this point of view, the modular invariants $ \sigma$ and $ \pi $ can be regarded as natural 
coordinates on the moduli space $\mathcal{M}_{{\rm M}}$.
\par Note then that the right-hand side space in $(\ref{hodgeident})$ also classifies unordered 
pairs $({\rm E}_1, {\rm E}_2) $ of curves of genus one. The two geometric structures, ${\rm M}$-polarized 
K3 surfaces and unordered pairs of elliptic curves, have the same classifying moduli space. Moreover, 
there is an obvious Hodge-theoretic bijective correspondence relating these structures:
\begin{equation}
\label{hodgecorresp}
\left ( {\rm X}, i \right )  \ \longleftrightarrow \ \left ( {\rm E}_1, {\rm E}_2 \right )
\end{equation} 
such that $ \sigma({\rm X}, i) = {\rm J}({\rm E}_1) + {\rm J}({\rm E}_2)$ and 
$\pi({\rm X}, i) = {\rm J}({\rm E}_1) \cdot {\rm J}({\rm E}_2)$.  

\section{A Geometric Transformation Underlying the Duality Map}
\label{partone}
As mentioned in the introduction, our goal is to place the the Hodge theoretic correspondence 
$ (\ref{hodgecorresp})$ into a geometric setting. In what follows, we provide the details 
needed for both the statement of Theorem $\ref{mainth}$ as well as its proof. 
Setting-up the geometric transformation requires a few technical ingredients 
concerning Shioda-Inose structures, Kummer surfaces and elliptic fibrations on a 
K3 surface. We shall therefore begin our exposition by presenting some basic facts. 

\subsection{Shioda-Inose Structures}
\label{shioda-inose}
The notion of a Shioda-Inose structure originates in the works \cite{shiodainose} of Shioda and 
Inose and \cite{nikulin} of Nikulin. Their ideas were later refined and generalized by Morrison 
\cite{morrison}. The above three papers are the main references for the assertions we review here.
\begin{dfn}
Let ${\rm X}$ be a K3 surface. An involution $\varphi \in {\rm Aut}({\rm X})$ is called a 
{\bf Nikulin involution} if $ \varphi^* \omega = \omega $ for any holomorphic two-form 
$\omega$. 
\end{dfn}
\noindent If a Nikulin involution $\varphi$ exists on ${\rm X}$, then $\varphi$ has exactly 
eight fixed points. In such a case, the quotient space 
$$ {\rm X} / \{ {\rm id}_{{\rm X}}, \varphi \} $$ 
is a surface with eight rational 
double point singularities of type $A_1$. The minimal resolution of this singular space is a new 
K3 surface which we denote by ${\rm Y}$. The two K3 surfaces ${\rm X}$ and ${\rm Y}$ are related 
by a (generically) two-to-one rational map $ \pi \colon {\rm X} \longrightarrow  {\rm Y}$. 
\par Denote by ${\rm H}^2_{{\rm Y}}$ the orthogonal complement in ${\rm H}^2({\rm Y}, \Zee) $ of the eight 
exceptional curves. One has then a natural push-forward map (see $\S$ 3 of \cite{morrison} or $\S$ 3 
of \cite{shiodainose}):
$$ \pi_* \colon {\rm H}^2({\rm X}, \Zee) \ \rightarrow \ {\rm H}^2_{{\rm Y}} $$ 
which restricts to a morphism of $\Zee$-modules:
\begin{equation}
\label{transcmorf}
\pi_* \colon {\rm T}_{{\rm X}} \ \rightarrow \ {\rm T}_{{\rm Y}} 
\end{equation}
between the transcendental lattices of the two K3 surfaces. 
\begin{rem}
\label{pairings}
The complexification of the morphism 
$ (\ref{transcmorf}) $ is a morphism of Hodge structures, but, in general, $ (\ref{transcmorf}) $ 
does not preserve the lattice pairings. In fact, one can check that:
$$ 
\langle \ \pi_* (t_1), \ \pi_*(t_2) \ \rangle _{{\rm Y}} \ = \ \langle \ t_1, t_2 \ \rangle _{{\rm X}} \ + \ 
\langle \ t_1, \ \varphi^*(t_2) \ \rangle _{{\rm X}}. 
$$
\end{rem} 
\begin{dfn}
A Nikulin involution $\varphi$ defines a {\bf Shioda-Inose structure} on ${\rm X}$ if ${\rm Y}$ 
is a Kummer surface and the morphism $(\ref{transcmorf})$ is a Hodge isometry 
${\rm T}_{{\rm X}}(2) \simeq {\rm T}_{{\rm Y}}$. 
\end{dfn} 
\noindent The notation ${\rm T}_{{\rm X}}(2)$ means that the bilinear pairing on the transcendental 
lattice ${\rm T}_{{\rm X}}$ is multiplied by $2$. We refer the reader to section $\ref{kummer}$ for 
an explanation of the significance of the last condition in the above definition, as well as for a 
short overview of the basics of Kummer surfaces.
\par Not every K3 surface admits a Nikulin involution, much less a Shioda-Inose structure. 
A very effective lattice-theoretic criterion which provides a necessary and sufficient condition 
for the existence of a Shioda-Inose structure on a K3 surface ${\rm X}$ has been given by Morrison. 
\begin{theorem} (Morrison \cite{morrison}, Theorem 5.7)
\label{morcr}
Let ${\rm X}$ be an algebraic K3 surface. There exists a Shioda-Inose structure on ${\rm X}$ if and only if the lattice 
$ {\rm E}_8 \oplus {\rm E}_8$ can be primitively embedded into the N\'{e}ron-Severi lattice 
${\rm NS}({\rm X})$.
\end{theorem} 
\noindent The proof of the above statement is based on a result of Nikulin (\cite{nikulin}, Theorem 4.3). 
A primitive embedding ${\rm E}_8 \oplus {\rm E}_8 \hookrightarrow {\rm NS}({\rm X})$ allows one to define 
a special lattice isometry of ${\rm H}^2({\rm X}, \Zee)$ which interchanges the two copies of 
${\rm E}_8$ given by the embedding and acts trivially on their orthogonal complement. In this context, 
Nikulin's theorem asserts that, possibly after conjugation by a reflection in an algebraic class of square $-2$, this lattice isometry can be associated to an involution of the K3 surface ${\rm X}$. Morrison shows then that this involution 
defines in fact a Shioda-Inose structure on ${\rm X}$.
\par Closer to the purpose of this paper, note that Theorem \ref{morcr} implies that 
an ${\rm M}$-polarized K3 surface $({\rm X},i)$ admits a Shioda-Inose structure. In fact, 
there exists a well-defined Shioda-Inose structure $\beta$ on ${\rm X}$ canonically associated 
with the ${\rm M}$-polarization. 
\par This canonical Shioda-Inose structure associated to an ${\rm M}$-polarized K3 surface 
plays a central role in our construction. However, in this paper we shall take a different point of view 
towards defining the Nikulin involution $\beta$ underlying the Shioda-Inose structure. Instead of using 
Theorem \ref{morcr}, we shall introduce this involution in a more explicit and geometric manner. The 
canonical involution $\beta$ appears naturally in the context of a special jacobian fibration on ${\rm X}$.    

\subsection{Jacobian Fibrations on K3 Surfaces}
\label{nikulin11}
During the course of this section we shall assume that ${\rm X}$ is an algebraic K3 surface.
\begin{dfn}
A {\bf jacobian fibration} (or elliptic fibration with section) on ${\rm X}$ is a pair $(\varphi, S)$ consisting of a proper map of analytic spaces 
$ \varphi \colon {\rm X} \rightarrow \mathbb{P}^1 $ whose generic fiber is a smooth curve of genus one, and a section 
$S$ in the elliptic fibration $\varphi$. 
\end{dfn}
\noindent If $S'$ is another section of the jacobian fibration $(\varphi, S)$, then there 
exists\footnote{See, for instance, Chapter 1 of \cite{morganfriedman} for a proof of this result.} an automorphism of ${\rm X}$ preserving $\varphi$ 
and mapping $S$ to $S'$. One can therefore realize an identification between the set of 
sections of $\varphi$ and the group of automorphisms of ${\rm X}$ preserving $\varphi$. 
This is the {\bf Mordell-Weil group} ${\rm MW}(\varphi, S)$ of the jacobian fibration.
\par Note also that a jacobian fibration $(\varphi, S)$ on ${\rm X}$ induces a sublattice:
$$ \mathcal{H}_{(\varphi,S)} \ \subset \ {\rm NS}({\rm X}) $$
constructed as the span of the two cohomology classes associated with the elliptic fiber and the 
section, respectively. The lattice $\mathcal{H}_{(\varphi,S)}$ is isomorphic to the standard rank-two 
hyperbolic lattice ${\rm H}$.
\par The sublattice $\mathcal{H}_{(\varphi,S)}$ determines uniquely the jacobian fibration $(\varphi,S)$. 
In other words, there cannot be two distinct jacobian fibrations on ${\rm X}$ 
determining the same hyperbolic sublattice in ${\rm NS}({\rm X})$. However, it is not true that any 
lattice embedding of ${\rm H}$ into ${\rm NS}({\rm X})$ corresponds to a jacobian fibration. Nevertheless, the following assertions hold: 
\begin{lem}
A lattice embedding ${\rm H} \hookrightarrow {\rm NS}({\rm X})$ can be associated with a jacobian fibration $(\varphi,S)$ if and only if its image in ${\rm NS}({\rm X})$ contains a pseudo-ample class.
\end{lem}
\begin{lem}
Let $\Gamma_{{\rm X}}$ be the group of isometries of ${\rm H}^2({\rm X}, \Zee)$ preserving the Hodge 
decomposition. For any lattice embedding 
$$ e \colon {\rm H} \hookrightarrow {\rm NS}({\rm X}),$$ 
there exists $ \alpha \in \Gamma_{{\rm X}} $ such that ${\rm Im}(\alpha \circ e) $ contains a pseudo-ample class.
\end{lem}
\begin{lem}
\label{l3}
One has the following bijective correspondence:  
\begin{equation}
\label{diag223}
\xymatrix{
\left \{  \txt{isomorphism classes of\\
jacobian fibrations on $X$ } \
\right \} 
\ \ar @{<->}  [r] & \    
\left \{ \
\txt{lattice embeddings\\ 
${\rm H} \hookrightarrow {\rm NS}(X)$}  \right \} / \ \Gamma_X
}. 
\end{equation}
\end{lem} 
\noindent These are standard well-known results. For proofs, we refer the reader to 
\cite{shapiro}, \cite{kondo} and \cite{clingher3}.
\par Next, let us consider 
$$ 
\mathcal{W}_{(\varphi, S)} \ \subset \ {\rm NS}({\rm X}) $$ 
to be the orthogonal complement of $\mathcal{H}_{(\varphi, S)}$ in the N\'{e}ron-Severi lattice of ${\rm X}$. 
It follows that $\mathcal{W}_{(\varphi, S)}$ itself is a negative-definite lattice 
of rank ${\rm p}_{{\rm X}}-2$. Moreover, 
the N\'{e}ron-Severi lattice decomposes as an orthogonal direct sum:
$$ {\rm NS}({\rm X}) \ = \ \mathcal{H}_{(\varphi, S)} \ \oplus \ \mathcal{W}_{(\varphi, S)}. $$
Let $ \Sigma \subset \mathbb{P}^1 $ be the set of points on the base of the elliptic fibration 
$ \varphi $ which correspond to singular fibers. For each $ v \in \Sigma $, denote by  ${\rm T}_v$ 
the sublattice of $\mathcal{W}_{(\varphi, S)}$ spanned by the classes of the irreducible components 
of the singular fiber over $v$ which are disjoint from $S$. One has then the following result. 
\begin{lem}
\label{jacoblem} \ 

\begin{itemize}
\item [(a)] For each $v \in \Sigma$, ${\rm T}_v$ is a negative-definite 
lattice of ${\rm ADE}$ type. 
\item [(b)] Let $ \mathcal{W}^{{\rm root}}_{(\varphi, S)} $ be the lattice spanned by the roots\footnote{A root of 
${\rm NS}({\rm X})$ is an algebraic class of self-intersection $-2$.} of 
$ \mathcal{W}_{(\varphi, S)} $. Then:
\begin{equation}
\label{decompk}
\mathcal{W}^{{\rm root}}_{(\varphi, S)} \ = \ \bigoplus _{v \in \Sigma} \ {\rm T}_v. 
\end{equation}
The decomposition $(\ref{decompk})$ is unique, up to a permutation of the factors.
\item [(c)] There exists a canonical group isomorphism:
\begin{equation}
\label{mwizo}
 \mathcal{W}_{(\varphi, S)} / \mathcal{W}^{{\rm root}}_{(\varphi, S)} \ 
 \stackrel{\simeq}{\longrightarrow} \ {\rm MW}(\varphi).
 \end{equation}
\end{itemize}
\end{lem}
\noindent The first two statements of the above lemma are standard facts from Kodaira's classification 
of singular fibers of elliptic fibrations (see, for example, \cite{kodaira1}). The last statement is due 
to Shioda \cite{shioda2}. 
\par Let us briefly indicate the construction of the correspondence 
in $(\ref{mwizo})$. Given 
$\gamma \in \mathcal{W}_{(\varphi, S)} $, 
denote by ${\rm L}$ the unique holomorphic line bundle over ${\rm X}$ such that 
$c_1({\rm L}) = \gamma $. Let $x \in {\rm X}$ be a point belonging to a smooth 
fiber ${\rm E}_{\varphi(x)}$. Then, the restriction of ${\rm L}$ to ${\rm E}_{\varphi(x)}$ is a holomorphic 
line bundle of degree zero and, therefore, there exists a unique $y \in E_{\varphi(x)} $ such that:
$$ {\rm L} \vert _{{\rm E}_{\varphi(x)}} \ \simeq \ \mathcal{O}_{{\rm E}_{\varphi(x)}}(x-y).$$
The assignment $x \mapsto y$ extends by continuity to an automorphism of the K3 surface and hence to an element 
in ${\rm MW}(\varphi)$.
\subsection{A Canonical Involution} 
We shall apply now the general theory presented in the previous section in the context of an ${\rm M}$-polarized 
K3 surface $({\rm X}, i)$.  
\par By standard lattice theory (see, for example, \cite{nikulin}), there are exactly two distinct 
ways (up to an overall isometry) in which one can embed the standard rank-two hyperbolic lattice 
${\rm H}$ isometrically into ${\rm M}$. The two possibilities are distinguished by the 
isomorphism type of the orthogonal complement of the image of the embedding. The orthogonal complement 
has rank $16$ and it is also unimodular, even, and negative-definite. As is well-known, up to isomorphism 
there exist only two such lattices. One is ${\rm E}_8 \oplus {\rm E}_8$ and the other is the Barnes-Wall lattice 
$ {\rm D}^+_{16} $. 
\par In the presence of an ${\rm M}$-polarization on ${\rm X}$,  
the two distinct isometric embeddings of the rank-two hyperbolic lattice ${\rm H}$ into ${\rm M}$ 
determine two distinct classes of embeddings of ${\rm H}$ into the N\'{e}ron-Severi lattice ${\rm NS}(X)$. 
According to Lemma $\ref{l3}$, one obtains therefore two special jacobian fibrations $(\Theta_1, S_1)$ and 
$(\Theta_2, S_2)$ on ${\rm X}$.
$$ \Theta_1, \Theta_2  \colon X \rightarrow \mathbb{P}^1. $$ 
We shall use the term {\bf standard fibration} for $\Theta_1$ (associated to 
the rank-sixteen lattice ${\rm E}_8 \oplus {\rm E}_8 $) and {\bf alternate fibration} for $\Theta_2$ 
(associated to the rank-sixteen lattice $ {\rm D}^+_{16}$). 
\begin{prop}
\label{prop11}
Let $({\rm X}, i)$ be an ${\rm M}$-polarized K3 surface.
\begin{itemize}
\item [(a)] The standard fibration $(\Theta_1, S_1) $ has two singular fibers of Kodaira type ${\rm II}^*$. The 
section $S_1$ is the unique section of $\Theta_1$ whose cohomology class belongs to $i({\rm M})$. 
\item [(b)] The alternate fibration $(\Theta_2, S_2) $ has a singular fiber of type ${\rm I}^*_{12}$. 
There are precisely two sections $S_2$ and $S'_2$ of $\Theta_2$ with cohomology classes represented 
in $i({\rm M})$. $S_2$ and $S'_2$ are disjoint. Moreover, the Mordell-Weil group $ {\rm MW}(\Theta_2)$ 
contains a canonical involution $\beta \in {\rm Aut}(X)$ which exchanges $S_2$ and $S'_2$.
\end{itemize}
\end{prop}
\bpf
The above assertions are consequences of the general principles reviewed in Section $\ref{nikulin11}$. 
In the case of the standard fibration $(\Theta_1, S_1)$, one has an orthogonal decomposition:
$$ \mathcal{W}^{{\rm root}}_{(\Theta_1, S_1)} \ = \ {\rm E}_8 \oplus {\rm E}_8 \oplus \mathcal{U}^{{\rm root}} $$
where $\mathcal{U}$ is the orthogonal complement of $i({\rm M})$ in ${\rm NS}({\rm X})$ and 
$\mathcal{U}^{{\rm root}}$ is the root lattice of $\mathcal{U}$. The above decomposition, combined with 
assertion (b) of Lemma $\ref{jacoblem}$, proves the existence of two singular fibers of Kodaira type 
${\rm II}^*$ in the elliptic fibration $\Theta_1$. It also follows immediately that $S_1$ is the unique section 
of $\Theta_1$ with associated cohomology class in $i({\rm M})$. One can represent the rational curves obtained as irreducible components of the 
two ${\rm II}^*$ fibers of $\Theta_1$ as well as the section $S_1$ in the following dual diagram.
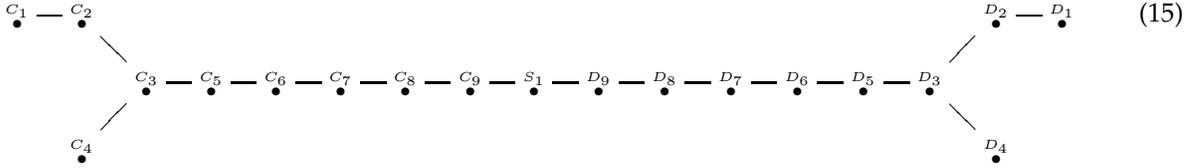
\begin{equation}
\label{diag11}
\def\objectstyle{\scriptstyle}
\def\labelstyle{\scriptstyle}
\xymatrix @-1.2pc  {
\stackrel{C_1}{\bullet} \ar @{-} [r] & \stackrel{C_2}{\bullet} \ar @{-} [dr] & & & & & & & & & & & & & & \stackrel{D_{2}}{\bullet} \ar @{-} [dl] 
 & \stackrel{D_{1}}{\bullet} \ar @{-} [l] \\
& & \stackrel{C_3}{\bullet} \ar @{-} [r] \ar @{-} [dl] &
\stackrel{C_5}{\bullet} \ar @{-} [r] &
\stackrel{C_6}{\bullet} \ar @{-} [r] &
\stackrel{C_7}{\bullet} \ar @{-} [r] &
\stackrel{C_8}{\bullet} \ar @{-} [r] &
\stackrel{C_9}{\bullet} \ar @{-} [r] &
\stackrel{S_1}{\bullet} \ar @{-} [r] &
\stackrel{D_{9}}{\bullet} \ar @{-} [r] &
\stackrel{D_{8}}{\bullet} \ar @{-} [r] &
\stackrel{D_{7}}{\bullet} \ar @{-} [r] &
\stackrel{D_{6}}{\bullet} \ar @{-} [r] &
\stackrel{D_{5}}{\bullet} \ar @{-} [r] &
\stackrel{D_{3}}{\bullet} \ar @{-} [dr] & \\
 & \stackrel{C_{4}}{\bullet} & & & & & & &  & & &  &  & & & \stackrel{D_{4}}{\bullet} \\
}
\end{equation}
\par The case of the alternate fibration $(\Theta_2, S_2)$ can be handled similarly. In this situation, 
one obtains an orthogonal decomposition:
$$  \mathcal{W}_{(\Theta_2, S_2)} \ = \ {\rm D}_{16}^+ \oplus \mathcal{U}. $$
Since any given root of the above lattice has to lie in one of the two factors, one obtains:
$$ \mathcal{W}^{{\rm root}}_{(\Theta_2, S)}  \ = \ {\rm D}_{16} \oplus \mathcal{U}^{\rm root}. $$
Once more, the assertion (b) of Lemma $\ref{jacoblem}$ tells one that $\Theta_2$ has a singular fiber 
of type ${\rm I}_{12}^*$. Next, note that, by assertion (c) of Lemma $\ref{jacoblem}$, 
one has an isomorphism of groups:
\begin{equation}
\label{mmm}
{\rm MW}(\Theta_2) \ \simeq \ \mathbb{Z}/ 2 \mathbb{Z} \ \oplus \  \mathcal{V}/ \mathcal{V}^{{\rm root}}. 
\end{equation}
The image $\beta \in {\rm MW}(\Theta_2)$ of the generator of the $\mathbb{Z}/ 2 \mathbb{Z} $ factor 
above determines naturally a non-trivial canonical involution of the K3 surface ${\rm X}$. In particular, 
the jacobian fibration $(\Theta_2, S_2)$ has an extra section $S_2'$, the image of $S_2$ through $\beta$. 
One can easily see from $(\ref{mmm})$ that $S_2$ and $S_2'$ are the only sections of the elliptic fibration 
$\Theta_2$ with cohomology classes represented in the polarizing lattice $i({\rm M})$.
\par In fact, one can clearly see the special ${\rm I}_{12}^*$ singular fiber together with two sections 
in the dual diagram $(\ref{diag11})$. This special singular fiber of $\Theta_2$ is given by the divisor:
$$ \ C_2 + C_4 + 2 \left ( C_3 + C_4 + \cdots C_9 + S + D_9 + D_8 + \cdots D_3  \right ) + D_{4} + D_{2}. $$
whereas the two sections $S_2$ and $S_2'$ are represented by the two extremal curves $C_1$ and $D_1$.  
\epf \ 

\begin{rem}
\label{diagremark}
The effect of the involution $\beta$ on the diagram $ (\ref{diag11})$ amounts to a 
right-left flip which sends the C-curves to the corresponding symmetric D-curves and vice-versa. 
In particular, the restriction of $\beta$ to the middle rational curve ${\rm S}$ is a non-trivial 
involution of ${\rm S}$ with two distinct fixed points.  
\end{rem}
\begin{rem}
\label{fixtrans}
The induced morphism $\beta^* \colon {\rm H}^2({\rm X}, \Zee) \rightarrow {\rm H}^2({\rm X}, \Zee)$ 
restricts to the identity on the orthogonal complement of $i({\rm M})$. In particular, $\beta^*$ acts 
trivially on the transcendental lattice ${\rm T}_{{\rm X}}$.
\end{rem}

\noindent One may guess now that it is the canonical involution $\beta$ of Proposition $\ref{prop11}$ that gives rise 
to the canonical Shioda-Inose structure on ${\rm X}$ we mentioned at the end of Section $\ref{shioda-inose}$.
\newpage
\par We are now in position to formulate the main result of the paper: 
\begin{theorem}
\label{th111}
Let $({\rm X}, i)$ be an ${\rm M}$-polarized K3 surface.
\begin{itemize}
\item [(a)] The involution $\beta$ introduced above defines a Shioda-Inose structure on ${\rm X}$. 
\item [(b)] The minimal resolution ${\rm Y}$ of the quotient ${\rm X} / \beta $ is a K3 surface 
with a canonical Kummer structure. This structure realizes ${\rm Y}$ as the Kummer surface of 
an abelian surface ${\rm A} = {\rm E}_1 \times {\rm E}_2 $ canonically represented as a cartesian 
product of two elliptic curves. The two elliptic curves are unique, up to permutation.  
\item [(d)] The above geometric transformation induces a canonical Hodge isomorphism 
between the ${\rm M}$-polarized Hodge structure of ${\rm X}$ and the natural ${\rm H}$-polarized 
Hodge structure of ${\rm A}$.
\end{itemize}  
\end{theorem}
\noindent Before embarking on the proof of the above theorem, let us comment briefly on the two 
special jacobian fibrations $\Theta_1$ and $\Theta_2$ which we have uncovered in this section.
These two jacobian fibrations\footnote{There is also an interesting toric reinterpretation of 
$\Theta_1$ and $\Theta_2$. They are induced from toric fibrations on a particular toric 
Fano three-fold by restriction to the anti-canonical hypersurface. These two toric fibrations 
are beautifully illustrated in Figure 1 of \cite{candelas}.
} are canonically associated to an ${\rm M}$-polarization on a 
K3 surface ${\rm X}$. However, so far, it is the standard fibration $\Theta_1$ that has received 
the lion's share of attention in the literature\footnote{This is also the reason why we decided to use the 
terms {\it standard} for $\Theta_1$ and {\it alternate} for $\Theta_2$.}. 
An analysis of $\Theta_1$ appears in the original work of 
Inose \cite{inose} and, over the last ten years, $\Theta_1$ has been extensively studied in the 
string theory literature due to its connection with the ${\rm E}_8 \oplus {\rm E}_8$ heterotic string 
theory in eight dimensions. The alternate fibration $\Theta_2$, however, has been largely overlooked. 
Nevertheless, it is $\Theta_2$, with its non-trivial Mordell-Weil group, that gives rise to 
a canonical Shioda-Inose structure on the ${\rm M}$-polarized K3 surface ${\rm X}$ and leads 
one to a geometric explanation for the Hodge-theoretic duality map $(\ref{hodgecorresp})$. 
The alternate fibration $\Theta_2$ will play a central role in the remainder of this paper.

\subsection{Kummer Surfaces}
\label{kummer}
In order to give a proof of Theorem $\ref{th111}$, we shall need a few classical results concerning 
the geometry of Kummer surfaces. For detailed proofs of the facts mentioned in this brief review 
we refer the reader to \cite{nikulin3}, \cite{shapiro} and \cite{morrison}.
\par Let ${\rm A}$ be a two-dimensional complex torus. Such a surface is naturally endowed with an abelian group 
structure. One can consider therefore on ${\rm A}$ the special involution given by $-{\rm id}$. The fixed 
locus of $-{\rm id}$ consists of sixteen distinct points. Therefore the quotient:
\begin{equation}
\label{kummersurface}
{\rm A}/ \{ \pm {\rm id} \} 
\end{equation}
is a singular surface with sixteen rational double point singularities of type ${\rm A}_1$. It is well-known
that the minimal resolution of $(\ref{kummersurface})$ is a special K3 surface ${\rm Km}({\rm A})$ called the 
{\bf Kummer surface} of ${\rm A}$. 
\par As a first important feature of Kummer surfaces, we note that the 
Hodge structures of ${\rm A}$ and ${\rm Km}({\rm A})$ 
are closely related. Indeed, denote by $ p \colon {\rm A} \rightarrow {\rm Km}({\rm A})$ 
the rational map induced by the quotienting and resolution procedure 
described above. Then, as explained for example in \cite{morrison}, one has a natural morphism 
$$p_* \colon {\rm H}^2({\rm A}, \Zee) \rightarrow {\rm H}^2_{{\rm Km}({\rm A})}$$ 
where 
${\rm H}^2_{{\rm Km}({\rm A})}$ is the sublattice of ${\rm H}^2({\rm Km}({\rm A}), \Zee)$ of classes 
orthogonal to all the sixteen exceptional curves. The complexification of $p_*$ sends the class
of a holomorphic two-form on ${\rm A}$ to a class representing a holomorphic two-form on ${\rm Km}({\rm A})$ 
and, as an immediate consequence of Proposition 3.2 in \cite{morrison}, one obtains:
\begin{prop} 
\label{partprop1}
The map $p_*$ is an isomorphism and it induces a canonical 
Hodge isometry 
\begin{equation}
\label{pistar}
{\rm H}^2({\rm A}, \Zee)(2) \ \stackrel{p_*}{\simeq}  \ {\rm H}^2_{{\rm Km}({\rm A})}.
\end{equation} 
Moreover, $p_*( {\rm T}_{{\rm A}} ) = {\rm T}_{{\rm Km}({\rm A})} $ and the above identification leads 
to a Hodge isometry at the level of transcendental lattices:
\begin{equation}
{\rm T}_{{\rm A}}(2) \ \stackrel{p_*}{\simeq} \ {\rm T}_{{\rm Km}({\rm A})}.
\end{equation}
\end{prop}
\noindent Kummer surfaces represent a large class of K3 surfaces. In fact, it is known 
(see for example \cite{shapiro}) that they form a dense subset in the moduli space of 
K3 surfaces. One would therefore like to have a criterion for determining whether 
a given K3 surface is Kummer. A very effective lattice-theoretic criterion for answering 
this question has been introduced by Nikulin in \cite{nikulin3}.
\begin{dfn} 
Let 
\begin{equation}
\label{skummer}
\mathcal{R} \ = \ \bigoplus _{i = 1} ^{16} \  \Zee x_i 
\end{equation}
be the rank-sixteen even lattice with bilinear form defined by $ (x_i,x_j) = -2 \delta_{ij} $. 
By definition, the {\bf Kummer lattice} $ {\rm K}$ is the lattice in $ \mathcal{R} \otimes \mathbb{Q} $ 
spanned (over $\Zee$) by:
$$ x_1, x_2, x_3, \cdots , x_{16} \ \ {\rm and} \ \ d \ = \ \frac{1}{2} \sum_{i=1}^{16}x_i, $$
\end{dfn}
\noindent The Kummer lattice ${\rm K}$ has rank sixteen, is even and negative-definite, and has 
the same discriminant group and discriminant form as the orthogonal sum:
$$ {\rm H}(2) \oplus {\rm H}(2) \oplus {\rm H}(2) $$
where ${\rm H}$ is the standard rank-two hyperbolic lattice. 
\par For any Kummer surface 
${\rm Km}({\rm A})$, one has a natural primitive lattice embedding:
$$ {\rm K} \ \hookrightarrow \ {\rm NS} \left ( {\rm Km}({\rm A}) \right ) $$
whose image is the minimal primitive sublattice of ${\rm NS} \left ( {\rm Km}({\rm A}) \right ) $ 
containing the classes of the sixteen exceptional curves. Nikulin's criterion asserts that 
the converse of the above statement is also true.
\begin{theorem}(Nikulin \cite{nikulin3}) 
\label{nikulincriterion}
\begin{itemize}
\item [(a)] A K3 surface ${\rm Y}$ is a Kummer surface if and only if there exists a primitive 
lattice embedding $ {\rm K} \hookrightarrow {\rm NS}({\rm Y}) $. 
\item [(b)] For every primitive lattice embedding $ e \colon {\rm K} \hookrightarrow {\rm NS}({\rm Y})$, 
there exists an unique and canonically defined two-dimensional complex torus ${\rm A}$ and a Hodge 
isometry $\alpha$ of ${\rm H}^2({\rm Y}, \Zee)$ such that ${\rm Y} = {\rm Km}({\rm A})$ 
and ${\rm Im}(\alpha \circ e)$ is the minimal primitive sublattice of ${\rm NS}({\rm Y})$ 
containing the sixteen exceptional curves arising during the Kummer construction process.
\end{itemize}
\end{theorem}
\noindent It is possible for a K3 surface ${\rm Y}$ to have multiple non-equivalent Kummer structures, i.e. there exist non-isomorphic complex tori ${\rm A}$ and 
${\rm A}'$ such that 
$$ {\rm Km}({\rm A}) \simeq {\rm Y} \simeq {\rm Km}({\rm A}').$$
However, as the last part of Theorem \ref{nikulincriterion} illustrates, once a primitive lattice embedding of the 
Kummer lattice ${\rm K}$ into ${\rm NS}({\rm Y})$ is fixed, there exists a unique complex torus ${\rm A}$ 
compatible with the embedding of ${\rm K}$. 
For a detailed treatment of the classification problem for Kummer structures 
on a K3 surface we refer the reader to the paper \cite{hosono3} of Hosono, Lian, Oguiso and Yau.
\par For the remainder of this section we shall restrict our attention to Kummer surfaces ${\rm Km}({\rm A})$ associated 
to abelian surfaces ${\rm A}={\rm E}_1 \times {\rm E}_2 $ realized as a cartesian product  
of two elliptic curves.  
\par Let us first introduce the basic criterion for an abelian surface ${\rm A}$ to have the above property. 
According to the Hodge index theorem, 
the N\'{e}ron-Severi lattice of ${\rm A}$, denoted by ${\rm NS}({\rm A})$, is an even lattice of 
signature $(1,r)$ with $0 \leq r \leq 3$. If ${\rm A}$ splits as a cartesian product 
${\rm E}_1 \times {\rm E}_2$ of two elliptic curves, then the cohomology classes of 
the two curves ${\rm E}_1$ and ${\rm E}_2$ span a rank-two hyperbolic 
sublattice of ${\rm NS}({\rm A})$. The converse of this statement also holds. 
\newpage
\begin{prop} 
\label{critsplit}
Let ${\rm A}$ be an abelian surface.
\begin{itemize}
\item [(a)] The surface ${\rm A}$ can be realized as a product of two elliptic curves if and only if 
there exists a primitive lattice embedding ${\rm H} \hookrightarrow {\rm NS}({\rm A})$.
\item [(b)] For every primitive lattice embedding $ e \colon {\rm H} \hookrightarrow {\rm NS}({\rm A})$, 
there exist two elliptic curves ${\rm E}_1$ and ${\rm E}_2$ (unique, up to permutation) and an 
analytic isomorphism ${\rm A} \simeq {\rm E}_1 \times {\rm E}_2$ such that ${\rm Im}(e)$ is spanned by 
the cohomology classes of ${\rm E}_1$ and ${\rm E}_2$.
\end{itemize}
\end{prop}  
\bpf 
This is a lattice-theoretic version of Ruppert's criterion for an abelian surface to be isomorphic 
to a cartesian product of two elliptic curves. See \cite{ruppert} or Chapter 10 $\S$ 6 of \cite{birkenhake} 
for proofs.
\epf \ 

\noindent Note that it is possible for an abelian surface ${\rm A}$ to be represented as a cartesian 
product of two elliptic curves in two or more non-equivalent ways. One can see that 
this phenomenon happens only when the Picard rank of ${\rm A}$ is maximal (${\rm p}_{{\rm A}}=4$). 
In such a case, the number of non-equivalent representations $ {\rm A} = {\rm E}_1 \times {\rm E}_2 $
has an interesting interpretation in the context of the class group theory of imaginary quadratic fields \cite{hosono2}. 
\par Let us assume now that a splitting ${\rm A} = {\rm E}_1 \times {\rm E}_2$ has been fixed. In this 
context, the cartesian product structure of ${\rm A}$ gives rise to a special configuration of 
twenty-four curves on the Kummer surface ${\rm Km}({\rm A})$. In order to introduce this curve configuration, 
let $\{ x_0, x_1,x_2, x_3\}$ and $\{ y_0,y_1,y_2,y_3\}$ be the two sets of points of order two on 
 ${\rm E}_1$ and ${\rm E}_2$. Denote by ${\rm H}_i$, ${\rm G}_j$ 
$(0 \leq i,j \leq 3)$ the rational curves on ${\rm Km}({\rm A})$ obtained as proper transforms 
of ${\rm E}_1 \times \{ y_i\} $ and $\{x_j\} \times {\rm E}_2$, respectively. Let also ${\rm E}_{ij}$ be the exceptional curve on 
${\rm Km}({\rm A})$ associated 
to the double point $(x_i,y_j)$ of ${\rm A}$. One has then the following intersection numbers:
$$ {\rm H}_i \cdot {\rm G}_j = 0$$ 
$$ {\rm H}_k \cdot {\rm E}_{ij} = \delta_{ki}, \ \  {\rm G}_k \cdot {\rm E}_{ij} = \delta_{kj}. $$
\begin{dfn}
\label{defdoublekummer}
The configuration of twenty-four rational curves 
\begin{equation}
\label{doublekummer1}
\{ {\rm H}_i, {\rm G}_j, {\rm E}_{ij} \ \vert \ 0 \leq i,j \leq 3\} 
\end{equation}
is called {\bf the double Kummer pencil} of ${\rm Km}({\rm A})$. The minimal primitive 
sublattice ${\rm NS}\left ( {\rm Km}({\rm A}) \right ) $ containing the classes 
of the curves in $(\ref{doublekummer1})$ is called {\bf the double Kummer lattice} of ${\rm Km}({\rm A})$. 
We denote the isomorphism class of this lattice by ${\rm DK}$. 
\end{dfn}
\begin{rem}
\label{uniquedk} By standard lattice theory, up to an overall isometry, 
the double Kummer lattice ${\rm DK}$ has a unique primitive embedding into the 
K3 lattice. The orthogonal complement of any such embedding is isomorphic to 
$${\rm H}(2) \oplus {\rm H}(2).$$
\end{rem}
\noindent Let us also note that the double Kummer lattice ${\rm DK}$ contains a natural 
finite-index sublattice isomorphic to 
$$ {\rm K} \oplus {\rm H}(2). $$ 
The left-hand side term above is, of course, the minimal primitive 
sublattice of ${\rm NS}\left ( {\rm Km}({\rm A}) \right ) $ 
containing the sixteen exceptional curves ${\rm E}_{ij}$ whereas the 
factor on the right-hand side is spanned by the two classes:
\begin{equation}
\label{2classes}
2{\rm H}_i + \sum_{j=0}^3 \ {\rm E}_{ij}, \ \ \ 2{\rm G}_j + \sum_{i=0}^3 \ {\rm E}_{ij}. 
\end{equation}
The two classes described above do not depend on the indices $i$ and $j$, respectively. Moreover, one can 
verify that the two classes of $(\ref{2classes})$ are precisely the images of the cohomology classes in 
${\rm H}^2({\rm A}, \Zee)$ associated to ${\rm E}_1$ and ${\rm E}_2$ under the morphism $ \pi_* $ of 
Proposition $\ref{partprop1}$.   
\par To summarize, we have seen that every Kummer surface $ {\rm Z}= {\rm Km}({\rm E}_1 \times {\rm E}_2)$ 
associated to an abelian surface that can be realized as a cartesian product of two genus-one curves 
comes equipped with a natural primitive lattice embedding ${\rm DK}  \hookrightarrow {\rm NS}({\rm Z})$.
In fact, one can see that the existence of such an embedding is a sufficient criterion for a K3 surface 
${\rm Z}$ to be a Kummer surface associated to a product abelian surface. 
\begin{prop}
\label{maintechpoint}
Let ${\rm Z}$ be a K3 surface. Assume that a primitive lattice embedding 
$e \colon {\rm DK} \hookrightarrow {\rm NS}({\rm Z})$ has been given. Then there exist two 
elliptic curves ${\rm E}_1 $ and ${\rm E}_2$ and a Hodge isometry $\alpha$ of ${\rm H}^2({\rm Z}, \Zee)$ 
such that 
$${\rm Z} = {\rm Km}({\rm E}_1 \times {\rm E}_2)$$ 
and ${\rm Im}(\alpha \circ e)$ is the 
double Kummer lattice associated to the Kummer construction. The two elliptic curves ${\rm E}_1$ 
and ${\rm E}_2$ are unique (up to permutation) and canonically defined. 
\end{prop} 
\bpf The above assertion is a consequence of the results presented earlier in this section. Let 
$$ e \colon {\rm DK} \hookrightarrow {\rm NS}({\rm Z})$$ 
be a primitive lattice embedding. By standard lattice theory, there exists a primitive 
embedding of the Kummer lattice ${\rm K}$ in $e({\rm DK})$. Moreover, this embedding is unique, 
up to an overall Hodge isometry of ${\rm H}^2({\rm Z}, \Zee)$. Therefore, by Nikulin's criterion, 
one has a canonical Kummer structure on ${\rm Z}$. In other words ${\rm Z}= {\rm Km}({\rm A})$ 
with ${\rm A}$ uniquely defined. Then, according to Proposition $\ref{partprop1}$, one has a
Hodge isometry 
\begin{equation}
\label{pistar22}
{\rm H}^2({\rm A}, \Zee)(2) \ \stackrel{p_*}{\simeq}  \ {\rm H}^2_{{\rm Z}}.
\end{equation}   
But ${\rm H}^2_{{\rm Z}}$ contains a canonical primitive sublattice of type ${\rm H}(2)$, the 
orthogonal complement of the Kummer lattice in $e({\rm DK})$. The preimage of this lattice 
in ${\rm H}^2({\rm A}, \Zee)$ is primitively embedded in ${\rm NS}({\rm A})$ and is isomorphic 
to ${\rm H}$. Then, by Proposition $\ref{critsplit}$, the abelian surface splits canonically as
a product of two elliptic curves. 
\epf   \

\subsection{Proof of Theorem \ref{th111}} 
\label{proofth1}
We are now in position to give detailed proofs for the statements of Theorem $\ref{th111}$. 
\par Let us begin by observing that $\beta$ is a Nikulin involution. If $\omega $ 
is a given holomorphic two-form on ${\rm X}$, then either 
$\beta^*\omega = \omega$ or $\beta^*\omega = - \omega$. But, the latter possibility implies 
(see, for example, \cite{zhang}) that either $\beta$ has no fixed locus (case that is ruled out 
by Remark $\ref{diagremark}$) or that the fixed locus of $\beta$ is a union of curves (case that is 
ruled out by the fact that $\beta$ acts without fixed points on the smooth fibers of $\Theta_2$). 
Therefore the only possibility that can occur is $\beta^*\omega=\omega$ which, by definition, 
means that $\beta$ is a Nikulin involution.   
\begin{rem}
As is well-known (for a proof of this fact see $\S 5$ of \cite{nikulin2}), the fixed locus of a Nikulin 
involution always consists of eight distinct points. The eight fixed points associated to 
$\beta$ appear nicely in the context of the alternate fibration $\Theta_2$. As noted in Remark 
$\ref{diagremark}$, two of them lie on the smooth rational curve $S$ (the middle curve 
of the dual diagram $(\ref{diag11})$, also the section of the standard fibration $\Theta_1$). The 
additional six fixed points lie on the singular fibers of $\Theta_2$. For instance, in the generic 
case, the alternate elliptic fibration $\Theta_2$ has, in addition to the ${\rm I}^*_{12}$ fiber, another 
six singular fibers of Kodaira type ${\rm I}_1$ (each consisting of a reduced rational 
curve with one node). The extra six fixed points of $\beta$ are precisely the nodes of those fibers. 
\end{rem} 
\noindent Let then ${\rm Y}$ be the K3 surface obtained as the minimal resolution of 
the quotient of ${\rm X}$ through $\beta$. We show now that ${\rm Y}$ is a Kummer surface. 
In order to carry out our argument, we denote by 
${\rm F}_1, {\rm F}_2 \cdots {\rm F}_8$ the eight exceptional curves arising after resolving 
the eight rational singularities. Assume that ${\rm F}_1 $ and ${\rm F}_2$ are associated 
to the two fixed points of $\beta$ that lie on the ${\rm I}_{12}^*$ fiber of $\Theta_2$. 
\par Recall that the alternate elliptic fibration $\Theta_2$ is left invariant by the involution 
$\beta$. Therefore, $\Theta_2$ induces a new elliptic fibration on ${\rm Y}$. We denote this 
fibration by $\Psi_2$.
 \begin{equation}
\label{diag4444}
\def\objectstyle{\scriptstyle}
\def\labelstyle{\scriptstyle}
\xymatrix 
{
{\rm X} \ar [dr] _{\Theta_2} \ar @{-->}[rr] ^{\pi} &  & {\rm Y} \ar [dl] ^{\Psi_2} \\
& \mathbb{P}^1 & \\ 
}
\end{equation}
It is then not hard to see that the ${\rm I}_{12}^*$ fiber of $\Theta_2$ becomes a singular fiber of 
Kodaira type ${\rm I}^*_6$ in the fibration $\Psi_2$. We represent its irreducible components in 
the dual diagram below.  
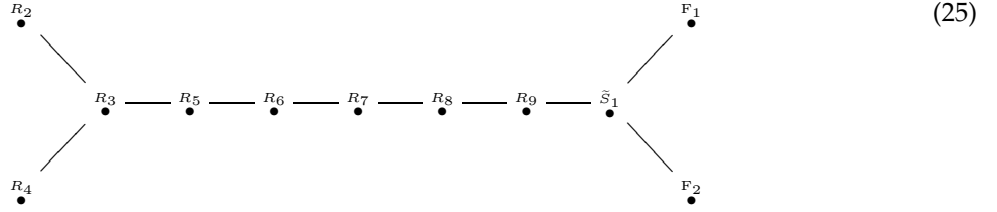
\begin{equation}
\label{diag22}
\def\objectstyle{\scriptstyle}
\def\labelstyle{\scriptstyle}
\xymatrix @-0.6pc  {
\stackrel{R_2}{\bullet} \ar @{-} [dr] & & & &  & &  & & \stackrel{{\rm F}_{1}}{\bullet} \ar @{-} [dl] 
 \\
 & \stackrel{R_3}{\bullet} \ar @{-} [r] \ar @{-} [dl] &
\stackrel{R_5}{\bullet} \ar @{-} [r] &
\stackrel{R_6}{\bullet} \ar @{-} [r] &
\stackrel{R_7}{\bullet} \ar @{-} [r] &
\stackrel{R_8}{\bullet} \ar @{-} [r] &
\stackrel{R_9}{\bullet} \ar @{-} [r] &
\stackrel{\w{S}_1}{\bullet}  \ar @{-} [dr] & \\
  \stackrel{R_{4}}{\bullet} & & & & & & & & \stackrel{{\rm F}_{2}}{\bullet} \\
}
\end{equation}
The curves ${\rm R}_i$, $1 \leq i \leq 9$ are the images of the curves 
${\rm C}_i $ (and also ${\rm D}_i$) of ${\rm X}$ (recall diagram $(\ref{diag11})$). 
The curve $\w{S}_1$ above is the quotient of the rational curve $S_1$ 
of diagram $(\ref{diag11})$ by the involution $\beta$. Note also that 
$R_1$ is a section in $\Psi_2$ while the unaccounted for exceptional 
curves ${\rm F}_3,{\rm F}_4 \cdots {\rm F}_8$ are disjoint from $R_1$ and form irreducible components 
in the additional singular fibers of $\Psi_2$.  
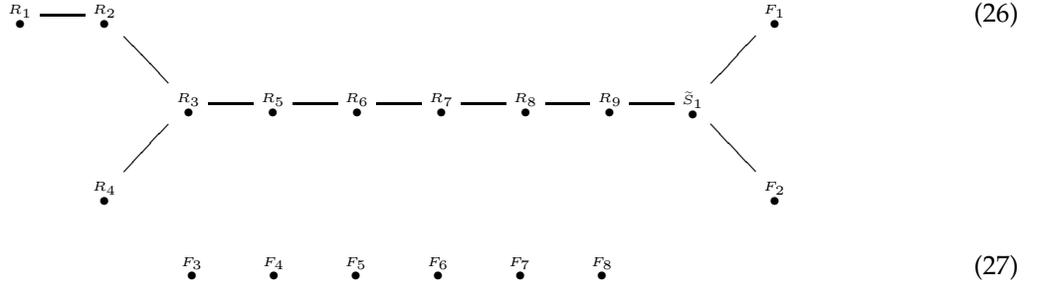
\begin{equation}
\label{diag33}
\def\objectstyle{\scriptstyle}
\def\labelstyle{\scriptstyle}
\xymatrix @-0.6pc  
{
\stackrel{R_1}{\bullet} \ar @{-} [r] & \stackrel{R_2}{\bullet} \ar @{-} [dr] & & & &  & &  & & \stackrel{F_{1}}{\bullet} \ar @{-} [dl] 
 \\
& & \stackrel{R_3}{\bullet} \ar @{-} [r] \ar @{-} [dl] &
\stackrel{R_5}{\bullet} \ar @{-} [r] &
\stackrel{R_6}{\bullet} \ar @{-} [r] &
\stackrel{R_7}{\bullet} \ar @{-} [r] &
\stackrel{R_8}{\bullet} \ar @{-} [r] &
\stackrel{R_9}{\bullet} \ar @{-} [r] &
\stackrel{\w{S}_1}{\bullet}  \ar @{-} [dr]   \\
&  \stackrel{R_{4}}{\bullet} & & & & & & & & \stackrel{F_{2}}{\bullet} 
}
\end{equation}
\begin{equation}
\label{3}
\def\objectstyle{\scriptstyle}
\def\labelstyle{\scriptstyle}
\xymatrix @-0.6pc  {
\stackrel{F_3}{\bullet}   &
\stackrel{F_4}{\bullet}  &
\stackrel{F_5}{\bullet}  &
\stackrel{F_6}{\bullet}  &
\stackrel{F_7}{\bullet}  &
\stackrel{F_8}{\bullet}  \\
}
\end{equation}

\begin{lem}
\label{lemma57}
Let $\mathcal{L}({\rm Y})$ be the minimal primitive sublattice of ${\rm NS}({\rm Y})$ containing 
the classes associated to the eighteen curves ${\rm R}_i$ ($1 \leq i \leq 9$), ${\rm F}_j$ ($1 \leq j \leq 8$) and 
$\w{S}_1$. The lattice $\mathcal{L}({\rm Y})$ is isomorphic to the double Kummer lattice ${\rm DK}$. 
\end{lem}
\bpf 
Let us denote by $\mathcal{N}$ the minimal primitive sublattice of ${\rm NS}({\rm Y})$ containing the eight 
exceptional curves $ {\rm F}_i $, $1 \leq i \leq 8$. The lattice $\mathcal{N}$ can also be regarded as
the span of the nine classes 
$$ {\rm F}_1, {\rm F}_2, \  \cdots \ {\rm F}_8, \ \ \frac{1}{2} \ \sum_{i=1}^8 \ {\rm F}_i .$$
This is the so-called {\bf Nikulin lattice} (see $\S 5$ of \cite{morrison}). It has rank eight and has the same discriminant group 
and discriminant form as ${\rm H}(2) \oplus {\rm H}(2) \oplus {\rm H}(2)$.
\par Denote by ${\rm H}^2_{{\rm Y}}$ the orthogonal complement of $\mathcal{N}$ in ${\rm H}^2({\rm Y}, \Zee)$.
Then, as we described in Section $\ref{shioda-inose}$, the Shioda-Inose construction induces a natural 
push-forward morphism:
$$ 
\pi_* \colon {\rm H}^2({\rm X}, \Zee) \rightarrow {\rm H}^2_{{\rm Y}} \hookrightarrow {\rm H}^2({\rm Y}, \Zee).
$$
Adopting the notation of diagram $(\ref{diag11})$, one has that
$ \pi_* (S_1) = 2 \w{S}_1+F_1+F_2$ and $ \pi_*(C_j)= \pi_*(D_j) = R_j$ for $1 \leq j \leq 9$. In particular, 
the image under $\pi_*$ of the class of the elliptic fiber of the standard fibration $\Theta_1$ on ${\rm X}$ 
is:
\begin{equation}
\label{fibertrans1}
2R_1+4R_2+6R_3+3R_4+5R_5+4R_6+3R_7+2R_8+R_9.
\end{equation}
We consider then the following sublattices of ${\rm H}_{{\rm Y}}^2$:
\begin{itemize}
\item [(a)] $\mathcal{E}$ is the span of the curves $ {\rm R}_j$, $1 \leq j \leq 8$.
\item [(b)] $\mathcal{H}$ is the span of $\pi_*(S_1) $ and $(\ref{fibertrans1})$.
\item [(c)] $\mathcal{Q} = \pi_* \left ( i({\rm M})^{\perp} \right )$. 
\end{itemize}
Using Remarks $\ref{pairings}$ and $\ref{fixtrans}$, we deduce that the three lattices above 
are orthogonal to each other. Moreover, $\mathcal{E}$ is isomorphic to ${\rm E}_8$ (hence unimodular), 
$\mathcal{H}$ is isomorphic to ${\rm H}(2)$ and $\mathcal{Q}$ is isomorphic to ${\rm H}(2) \oplus {\rm H}(2)$.
Hence, the discriminant of $ \mathcal{E} \oplus \mathcal{H} \oplus \mathcal{Q} $ is $2^6$. But 
the lattice ${\rm H}_{{\rm Y}}^2 $ has the same discriminant as its orthogonal complement $\mathcal{N}$ which, 
in turn, has discriminant $2^6$. Since clearly $ \mathcal{E} \oplus \mathcal{H} \oplus \mathcal{Q} $ is 
a sublattice of ${\rm H}_{{\rm Y}}^2 $, the equality of the two discriminants allows us to conclude that 
\begin{equation}
{\rm H}_{{\rm Y}}^2  \ = \ \mathcal{E} \oplus \mathcal{H} \oplus \mathcal{Q}.
\end{equation}
In particular, $\mathcal{Q}$ must be primitively embedded in ${\rm H}^2({\rm Y}, \Zee)$.
\par Now, by standard lattice theory (\cite{nikulin}, Theorem 1.14.4), up to an overall isometry there 
exists a unique primitive lattice embedding of ${\rm H}(2) \oplus {\rm H}(2)$ into the K3 lattice. 
By Remark $\ref{uniquedk}$, the orthogonal complement of such an embedding is isomorphic to the double Kummer 
lattice ${\rm DK}$. We see therefore that $\mathcal{Q}^{\perp} $ is isomorphic to ${\rm DK}$. 
\par At this point, let us also note the primitive embedding 
$ \mathcal{L} ({\rm Y}) \subset \mathcal{Q}^{\perp}$. In order to show that 
$ \mathcal{L} ({\rm Y}) = \mathcal{Q}^{\perp}$ all we need to do is verify that 
the two lattices involved have the same rank. The rank of $\mathcal{Q}^{\perp}$ is $18$, as it is 
isomorphic to ${\rm DK}$. By definition, rank$(\mathcal{L}({\rm Y})) \leq 18 $. But
$$ \mathcal{N} \oplus \mathcal{E} \oplus \mathcal{H} \ \subset \ \mathcal{L}({\rm Y}) $$ 
and therefore rank$(\mathcal{L}({\rm Y})) \geq 18 $. Hence, we have that 
$\mathcal{L}({\rm Y})= {\rm Q}^{\perp}$ and therefore $\mathcal{L}({\rm Y})$ 
is isomorphic with the double Kummer lattice ${\rm DK}$.
\epf \ 

\noindent The above result shows that, by construction, ${\rm Y}$ comes endowed with 
a canonical primitive lattice embedding ${\rm DK} \hookrightarrow {\rm H}^2({\rm Y}, \Zee)$. 
This fact, in connection with Proposition $\ref{maintechpoint}$, implies that there exist 
two canonically defined elliptic curves ${\rm E}_1$, ${\rm E}_2$ (unique, up to permutation) 
such that
\begin{equation}
\label{kummerc}
{\rm Y} \ = \ {\rm Km}({\rm E}_1 \times {\rm E}_2). 
\end{equation}
Moreover, Proposition $\ref{maintechpoint}$ together with weak form of the Global Torelli Theorem 
(Theorem 11.1 of \cite{bpv} $\S~VIII$), implies that $\mathcal{L}({\rm Y})$ is precisely the double Kummer 
lattice associated to the Kummer construction $(\ref{kummerc})$.
\par In order to check the last assumption of Theorem $\ref{th111}$, let us consider the diagram of rational 
maps:
\begin{equation}
\label{diag7878}
{\rm X} \ \stackrel{\pi}{\longrightarrow} \ {\rm Y} \ \stackrel{p}{\longleftarrow} \ {\rm E}_1 \times {\rm E}_2
\end{equation}
where $\pi$ is the map induced by the Shioda-Inose construction and $p$ is the map associated with 
the Kummer construction. The K3 surface ${\rm X}$ carries the lattice polarization 
$i({\rm M}) \subset {\rm H}^2({\rm X}, \Zee)$, whereas the abelian surface ${\rm E}_1 \times {\rm E}_2$ 
is ${\rm H}$-polarized by the sublattice ${\rm P} \subset {\rm H}^2({\rm E}_1 \times {\rm E}_2, \Zee)$ spanned by the 
classes of ${\rm E}_1$ and ${\rm E}_2$. In both cases, the orthogonal complement of the polarizing 
lattice is isomorphic to ${\rm H} \oplus {\rm H}$. One has then the push-forward morphisms:
\begin{equation}
{\rm H}^2({\rm X}, \Zee) \ \stackrel{\pi_*}{\longrightarrow} \ {\rm H}^2({\rm Y}, \Zee) \ 
\stackrel{p_*}{\longleftarrow} \ {\rm H}^2({\rm E}_1 \times {\rm E}_2, \Zee). 
\end{equation}
From the proof of Lemma $\ref{lemma57}$, we have that 
$\pi_*(i({\rm M})^{\perp}) = \mathcal{Q} = \mathcal{L}({\rm Y})^{\perp} $. By Proposition $\ref{partprop1}$, 
$ p_*({\rm P}^{\perp}) = \mathcal{L}({\rm Y})^{\perp} $. Moreover, both the restriction of the 
$\pi_*$ on $i({\rm M})^{\perp}$ and the restriction of $p_*$ on ${\rm P}^{\perp}$ induce isomorphisms 
of Hodge structures:
$$ 
i({\rm M})^{\perp}(2) \ \stackrel{\pi_*}{\simeq} \ \mathcal{L}({\rm Y})^{\perp}, \ \ \ \ 
{\rm P}^{\perp} (2) \ \stackrel{p_*}{\simeq} \ \mathcal{L}({\rm Y})^{\perp}. 
$$
By taking $(p_*)^{-1} \circ \pi_*$, one obtains therefore a canonical isomorphism of polarized Hodge structures
$$ i({\rm M})^{\perp} \ \stackrel{\pi_*}{\simeq} \ {\rm P}^{\perp}  $$
between the surfaces ${\rm X}$ and ${\rm E}_1 \times {\rm E}_2$.

\section{An Explicit Computation} 
\label{parttwo}
In the first half of this paper, we have described a geometric correspondence:
$$ ({\rm X}, i) \ \mapsto \ {\rm A}({\rm X}) = {\rm E}_1 \times {\rm E}_2 $$
which associates to any given ${\rm M}$-polarized K3 surface ${\rm X}$ an 
abelian surface ${\rm A}({\rm X})$ 
realized as a cartesian product of two elliptic curves. In this second part of the paper
we shall make this correspondence explicit. In other words we shall compute the J-invariants 
of the two elliptic curves ${\rm E}_1 $ and ${\rm E}_2$. 
\par Note that, by the Hodge theoretic equivalence underlying the correspondence, the modular 
invariants of an ${\rm M}$-polarized K3 surface $({\rm X},i)$ can be written as:
$$ \sigma({\rm X}, i) = {\rm J}({\rm E}_1) + {\rm J}({\rm E}_2), \ \ \ 
\pi({\rm X},i) = {\rm J}({\rm E}_1) \cdot {\rm J}({\rm E}_2). $$
Therefore, as an immediate application of the calculation of the two J-invariants of ${\rm E}_1$ and 
${\rm E}_2$, we shall obtain formulas 
for the modular invariants associated to an explicitly defined ${\rm M}$-polarized K3 surface.
\subsection{The Inose Form} 
\par In his 1977 paper \cite{inose}, Inose introduced an explicit two-parameter family of 
K3 surfaces which carry canonical ${\rm M}$-polarizations. The surfaces in this family are 
defined as follows. 
\par Let $a,b \in \Cee $. Denote by ${\rm Q}(a,b) $ the surface in $\mathbb{P}^3$ defined by the quartic equation:   
\begin{equation}
\label{inoseform1}
y^2zw - 4x^3z + 3axzw^2 - \frac{1}{2} \left ( z^2w^2 + w^4 \right )  + b zw^3 \ = \ 0 \ .
\end{equation}
We shall refer to the polynomial on the left side of the above equation as the {\bf Inose form}. 
The surface ${\rm Q}(a,b)$ has only rational double point singularities and its minimal resolution, denoted 
${\rm X}(a,b)$, is a K3 surface. Moreover, by construction, the surface ${\rm X}(a,b)$ has a canonical ${\rm M}$-polarization. In 
order to see this, let us note that the intersection of ${\rm Q}(a,b)$ with the hyperplane $\{ w=0\}$ 
is a union of two lines ${\rm L}_1 \cup {\rm L}_2$ with: 
$${\rm L}_1 : = \{ z=w=0\}, \ \ \ \ {\rm L}_2 := \{x=w=0\}. $$ 
Moreover, by standard singularity theory, the points $[0,1,0,0]$ and 
$[0,0,1,0]$ are rational double point singularities on ${\rm Q}(a,b)$ of types ${\rm A}_{11}$ and ${\rm E}_6$, 
respectively.  As a result, one obtains on the minimal resolution ${\rm X}(a,b)$ the 
following configuration of rational curves: 
\begin{equation}
\label{diaginose}
\def\objectstyle{\scriptstyle}
\def\labelstyle{\scriptstyle}
\xymatrix @-0.8pc  {
\stackrel{a_1}{\bullet} \ar @{-} [r]
& \stackrel{a_2}{\bullet} \ar @{-} [r]& \stackrel{a_3}{\bullet} \ar @{-} [r] \ar @{-} [d] &
\stackrel{a_4}{\bullet} \ar @{-} [r] &
\stackrel{a_5}{\bullet} \ar @{-} [r] &
\stackrel{a_6}{\bullet} \ar @{-} [r] &
\stackrel{a_7}{\bullet} \ar @{-} [r] &
\stackrel{a_8}{\bullet} \ar @{-} [r] &
\stackrel{a_9}{\bullet} \ar @{-} [r] &
\stackrel{a_{10}}{\bullet} \ar @{-} [r] &
\stackrel{a_{11}}{\bullet} \ar @{-} [r] &
\stackrel{L_{2}}{\bullet} \ar @{-} [r] &
\stackrel{e_{1}}{\bullet} \ar @{-} [r] &
\stackrel{e_{2}}{\bullet} \ar @{-} [r] &
\stackrel{e_{3}}{\bullet} \ar @{-} [d] & \stackrel{e_{5}}{\bullet} \ar @{-} [l] 
 & \stackrel{e_{6}}{\bullet} \ar @{-} [l] \\
 & & \stackrel{L_{1}}{\bullet} &  & & & & &  & & &    & & & \stackrel{e_{4}}{\bullet} & \\
}
\end{equation}
Note already the similarity with the previously encountered diagram $(\ref{diag11})$. The lattices 
spanned by:
$$ \{ \ a_1, \ a_2, \ L_1, \ a_3, \ a_4, \ a_5, \ a_6, \ a_7 \ \}$$
$$ \{ \ a_{11}, \ L_2, \ e_1, \ e_2, \ e_3, \ e_4, \ e_5, \ e_6 \ \}$$
$$ \{ \ a_9, \ 2a_1+4a_2+ 3L_1 + 6a_3 + 5 a_4 + 4 a_5 + 3 a_6 + 2 a_7 \ \} $$
are mutually orthogonal and they are also isomorphic to ${\rm E}_8$, ${\rm E}_8$ and ${\rm H}$, respectively. 
As a consequence, their direct sum provides a canonical primitive lattice embedding 
${\rm M} \hookrightarrow {\rm NS}({\rm X}(a,b))$.
\subsection{The Main Formula}
In the remaining part of the paper, we prove:
\begin{theorem}
\label{main}
Let ${\rm E}_1$ and ${\rm E}_2$ be the two elliptic curves associated to the ${\rm M}$-polarized K3 surface 
${\rm X}(a,b)$ by the correspondence of Theorem $\ref{th111}$. Then ${\rm J}({\rm E}_1)$ and ${\rm J}({\rm E}_2)$ 
are the two solutions of the quadratic equation:
\begin{equation}
\label{quadraticeq}
x^2 - \left ( a^3 -b^2 + 1\right )x + a^3  \ = \ 0. 
\end{equation} 
\end{theorem}
\noindent As mentioned earlier, as a consequence of the above theorem, one obtains:
\begin{cor}
The modular invariants of the ${\rm M}$-polarized K3 surface ${\rm X}(a,b)$ are given by:
\begin{equation}
\label{mainformula}
\pi = a^3 , \ \ \sigma  = a^3 - b^2 + 1.
\end{equation}
\end{cor}
\begin{cor}
Every ${\rm M}$-polarized K3 surface is isomorphic\footnote{Here by the term isomorphism we mean 
an isomorphism of ${\rm M}$-polarized K3 surfaces.} to ${\rm X}(a,b)$ for some $a,b \in \Cee$.
\end{cor}
\noindent Our strategy for proving Theorem $\ref{main}$ relies on a detailed analysis of the two 
basic algebraic invariants associated with the elliptic fibration $\Psi_2$ on the Kummer surface ${\rm Y}$:
the functional and homological invariants. 
\subsection{Invariants Associated to an Elliptic Surface}
\noindent Let $ {\rm X} $ be a smooth compact complex analytic surface and let 
$ \varphi \colon {\rm X} \rightarrow {\rm C}$ be a proper analytic map to a smooth curve 
such that the generic fiber of $\varphi$ is a smooth elliptic curve. Assume also that 
$\varphi$ does not have multiple fibers. There are two classical invariants that one 
associates to such a structure \cite{morganfriedman, kodaira1}.
\begin{itemize}
\item [(a)] The {\bf functional invariant} is an analytic function 
$\mathcal{J}_{\varphi} \colon {\rm C} \rightarrow \mathbb{P}^1$. It can be defined 
in the following manner. Let $ {\rm U}$ be the complement in ${\rm C}$ of the critical 
values of $\varphi$. Then $\mathcal{J}_{\varphi}$ is the meromorphic continuation of the 
composite map:
$$ {\rm U} \ \stackrel{e}{\longrightarrow} \ \mathbb{H}/{\rm PSL}(2, \Zee) \ \stackrel{{\rm J}}{\longrightarrow} \ \Cee $$
which takes a smooth elliptic fiber to its associated point in the moduli space of 
elliptic curves and then evaluates the classical elliptic modular function\footnote{Recall 
that ${\rm J}$ is normalized 
such that the two orbifold points of $\mathbb{H}/{\rm PSL}(2, \Zee)$ are mapped to $0$ and $1$.} at that respective point.
\item [(b)] The {\bf homological invariant} is, by definition, the sheaf 
$\mathcal{G}_{\varphi} = {\rm R}^1 \varphi _* \Zee_{{\rm X}} $. The restriction of $\mathcal{G}_{\varphi} $ on 
${\rm U}$ is locally constant and oriented and its stalk at every point is isomorphic with $\Zee \oplus \Zee$. 
Moreover, since $\varphi$ has no multiple fibers, one has $\mathcal{G}_{\varphi} = i_* \left ( \mathcal{G}_{\varphi}\vert _{{\rm U}} \right )  $ where $i \colon {\rm U} \hookrightarrow {\rm C} $ and 
therefore $\mathcal{G}_{\varphi}$ is determined by its restriction on ${\rm U}$. The latter sheaf is 
however fully determined by the conjugacy class of its monodromy map:
\begin{equation}
\label{monodromymap}
\rho_{\varphi} \colon \pi_1 \left ( {\rm U}, t \right ) \ \longrightarrow \ 
{\rm SO} \left ( {\rm H}^1(\varphi^{-1}(t), \Zee) \right )  
\end{equation}    
One can regard, therefore, the homological invariant of $\varphi$ as an element in 
${\rm Hom}\left ( \pi_1({\rm U}), {\rm SL}(2, \Zee) \right ) $, modulo conjugation. 
\end{itemize}  
\noindent The two invariants are not unrelated. Let us assume, for simplicity,  
that $\mathcal{J}_{\varphi}$ is not constant, as the cases of interest to us will definitely satisfy 
this condition. Set then ${\rm U}_0 \subset {\rm U}$ as the open subset for which 
$\mathcal{J}_{\varphi} \notin \{ 0,1\}$ and denote by $\mathbb{H}_0$ the set of elements of 
the upper half-plane $\mathbb{H}$ for which the associated elliptic modular function is neither 
$0$ or $1$. Pick $t \in {\rm U}_0$ . The composition 
$$ {\rm U}_0 \ \stackrel{i}{\hookrightarrow} \ {\rm U} \ \stackrel{e}{\rightarrow} \ \mathbb{H}/{\rm PSL}(2, \Zee)$$ 
induces a morphism of fundamental groups:
\begin{equation}
\label{morphism1}
 \pi_1\left ( {\rm U}_0, t \right ) \ \stackrel{(e \circ i)_{\#}}{\rightarrow} \ 
\pi_1 \left ( \mathbb{H}_0/{\rm PSL}(2, \Zee) \right ) \ \simeq \ {\rm PSL}(2, \Zee)
\end{equation} 
The compatibility between the two invariants asserts that the above morphism agrees, modulo conjugation with:
\begin{equation}
\label{morphism2}
 \pi_1\left ( {\rm U}_0, t \right ) \ \stackrel{i_{\#}}{\rightarrow} \  
 \pi_1\left ( {\rm U}, t \right ) \ \stackrel{\rho_{\varphi}}{\rightarrow} \ {\rm SL}(2, \Zee) \ \rightarrow {\rm PSL}(2, \Zee). 
\end{equation}
\par The above compatibility condition can be introduced independent of the actual elliptic 
fibration over ${\rm C}$. Given a (non-constant) meromorphic 
function $ \mathcal{J} \colon {\rm C} \rightarrow \Cee $ with no poles on ${\rm U}$ and a 
morphism of $\Zee$-modules  
$ \rho \colon \pi_1 \left ( {\rm U} \right ) \rightarrow {\rm SL}(2, \Zee)$, the pair $(\mathcal{J}, \rho)$ 
is said to be {\bf compatible} if the associated maps $(\ref{morphism1})$ and $(\ref{morphism2})$ agree modulo conjugation.  One has then the following classical theorem of Kodaira.
\begin{theorem} (Kodaira \cite{kodaira1})
For a compatible pair $(\mathcal{J},\rho)$ as above, there is, up to an isomorphism of elliptic surfaces, 
exactly one elliptic fibration $ \varphi \colon {\rm X} \rightarrow {\rm C}$, admitting a section, with functional and 
homological invariants given by $\mathcal{J}$ and $\rho$. 
\end{theorem}  
\noindent The above theorem provides one with a very powerful tool for comparing jacobian fibrations. Its 
effectiveness is further enhanced by the fact that, given a jacobian fibration as above, 
the monodromy $\rho(\gamma)$ of a small loop $\gamma$ circling a critical value of $\varphi$ in a 
manner agreeing with the orientation of ${\rm C}$ is determined modulo conjugation by 
the Kodaira type of the associated singular fiber \cite{kodaira2}. We can state therefore the following 
very particular consequence of the above discussion.
\begin{cor}
\label{cormain}
Let $ \varphi $ and $ \psi $ be two jacobian fibrations on two K3 surfaces ${\rm X}$ and ${\rm X'}$. 
 \begin{equation}
\label{chartt}
\def\objectstyle{\scriptstyle}
\def\labelstyle{\scriptstyle}
\xymatrix 
{
{\rm X} \ar [dr] _{\varphi}  &  & {\rm X'} \ar [dl] ^{\psi} \\
& \mathbb{P}^1 & \\ 
}
\end{equation}
The two jacobian fibrations are isomorphic if and only if there exists a a projective automorphism 
${\rm q}$ of $\mathbb{P}^1$ such that ${\rm q}$ maps bijectively the singular locus of $\varphi$ 
to the singular locus of $\psi$, $\mathcal{J}_{\varphi} = \mathcal{J}_{\psi} \circ {\rm q}$ and, for any 
$t$ in the singular locus of $\varphi$, the 
Kodaira type of a singular fiber $ \varphi^{-1}(t)$ is the same as the Kodaira type of 
the singular fiber $ \psi^{-1}({\rm q}(t))$. 
\end{cor}    
\noindent Our strategy for proving Theorem $\ref{main}$ is structured as follows. We first compute 
the functional invariants and Kodaira types of singular fibers of both the alternate fibration 
$\Theta_2$ on $ {\rm X}(a,b)$, and the induced jacobian fibration $ \Psi_2$ on the K3 surface 
${\rm Y}(a,b)$. Then, switching our attention to the other side of the correspondence, we 
show that, for any two elliptic curves ${\rm E}_1$ and ${\rm E}_2$, the Kummer surface 
${\rm Km}({\rm E}_1 \times {\rm E}_2)$ possesses a canonical jacobian fibration $\Upsilon_2$ 
with the same types of singular fibers as $\Psi_2$. Finally, using Corollary $\ref{cormain}$, 
we prove that the two elliptic fibrations $\Psi_2$ and $\Upsilon_2$ are equivalent if 
and only if ${\rm J}({\rm E}_1)$ and ${\rm J}({\rm E}_2)$ are solutions to the quadratic 
equation $(\ref{quadraticeq})$.
\subsection{The Alternate Fibration $\Theta_2$}
\label{step1}
It is quite easy to observe the fibration $\Theta_2$ on the surface ${\rm X}(a,b)$. The 
alternate fibration is induced by the projection to $[x,w]$ from the quartic ${\rm Q}(a,b)$.  Indeed, one can easily 
verify the following facts.
\begin{itemize}
\item [(a)] The generic fiber of the projection to $[x,w]$ from ${\rm Q}(a,b)$ is an elliptic curve. 
In fact, the fiber over $[1, \lambda]$ can be seen as the cubic curve in $\mathbb{P}^2(y,z,w)$ given by:
\begin{equation}
\label{cubics1}
\Theta_2^{\lambda} : = \Theta_2 ^{-1} \left ( [\lambda,1] \right ) \ = \ 
\{ \ 2y^2z - \left ( 8\lambda^3 - 6azw^2 - 2b \right ) z w^2 - z^2w -w^3  = 0 \}. 
\end{equation}  
This is a smooth cubic as long as  $4\lambda^3 -3a\lambda - b \neq \pm 1$.
\item [(b)] After resolving the singularities of ${\rm Q}(a,b)$, the projection to $[x,w]$ induces 
an elliptic fibration on the K3 surface ${\rm X}(A,B) $. 
\item [(c)] The singular fiber $\Theta_2^{\infty} : = \Theta_2^{-1} \left (  [1,0] \right )$ is of 
Kodaira type ${\rm I}_{12}^*$. In the context of the diagram $(\ref{diaginose})$, $\Theta_2^{\infty}$ 
appears as the divisor:
$$ a_2 + L_1 + 2 \left ( a_3 + a_4 + \cdots + a_{11} + L_2 + e_1 + e_2 + e_3 \right ) + e_4 + e_5 . $$
\item [(d)] The curves $a_1$ and $e_6$ are sections of $\Theta_2$.
\end{itemize} 
\noindent Let us then compute the functional invariant of the elliptic fibration $\Theta_2$. In order 
to simplify further calculations, we shall introduce the following polynomial: 
$$ {\rm P}(X) = 4X^3 - 3A X - B. $$ 
With this in place, one can rewrite the cubic equation in $(\ref{cubics1})$ in a 
standard Weierstrass form as:
\begin{equation}
\label{weierstrassf}
\left ( \sqrt{2} \ yz \right )^2 \ = \ \left ( z + \frac{2}{3} \ {\rm P}(\lambda) \right )^3 \ + \ 
 g_2(\lambda) \left ( z + \frac{2}{3} \ {\rm P}(\lambda) \right ) \ + \ g_3(\lambda) 
 \end{equation}
where the terms $g_2(\lambda)$ and $g_2(\lambda)$ are given by:
 $$ g_2(\lambda) = 1 - \frac{4}{3} \ {\rm P}^2(\lambda), \ \ \ 
 g_3(\lambda) = \frac{16}{27}\ {\rm P}^3(\lambda) - \frac{2}{3}\ {\rm P}(\lambda).  $$
The discriminant of Weierstrass form $(\ref{weierstrassf})$ is then:
$$ \Delta_{\Theta_2}(\lambda) \ = \ 
4 g_2^3(\lambda) + 27 g_3^2(\lambda) \ = \ 4 \left ( 1 - {\rm P}^2(\lambda) \right ). $$
In the same manner, the functional invariant of $ \Theta_2$ can be computed as:
$$ \mathcal{J}_{\Theta_2}(\lambda) \ = \ 
\frac{4g^3_2(\lambda)}{\Delta_{\Theta_2}(\lambda)}
\ = \ 
\frac{\left ( 3 - 4 \ {\rm P}^2(\lambda) \right )^2}{9(1 - {\rm P}^2(\lambda)) }. $$
The explicit formulas for $g_2(\lambda)$, $g_3(\lambda)$ and $\Delta_{\Theta_2}(\lambda)$ allow 
one to determine not only the location but also the Kodaira type of the singular fibers of $\Theta_2$. 
Using Tate's algorithm \cite{tate}, one obtains:
\begin{prop} 
\label{homol}
The singular fibers of $\Theta_2$ are located at $[1,0]$ (the ${\rm I}_{12}^*$ fiber) and at the 
points $[\lambda,1]$ with $\lambda$ belonging to the subset:
\begin{equation}
\label{singloc}
\Sigma := \ \{ \ \lambda \  \vert \ {\rm P}(\lambda)^2 = 1 \ \}. 
\end{equation}
The following five cases can occur:
\begin{itemize}
\item $a^3 \neq (b \pm 1)^2 $. In this case, both polynomials ${\rm P}(X)-1$ and ${\rm P}(X)+1$ have 
three distinct roots. The subset $\Sigma$ consists of six distinct points, each of which corresponds 
to a singular fiber of type ${\rm I}_1$ in $\Theta_2$.
\item $a^3 = (b + 1)^2 $,  $ b \neq 0$, $ a \neq 0$. In this case, ${\rm P}(X)+1$ has three distinct roots. 
However, the polynomial ${\rm P}(X)-1$ has a root of order two at $-(b+1)/2a$ and a simple root at $(b+1)/a$. 
The subset $\Sigma$ consists of 5 distinct points. 
$$ \Sigma \ = \ \left \{ \frac{-(b+1)}{2a}, \ \frac{b+1}{a} \right \} \ \cup \ 
\left \{ \ \lambda \ \vert \ {\rm P}(\lambda)=-1 \ \right \} $$
The first value in the above list corresponds to a singular fiber of type 
${\rm I}_2$ in $\Theta_2$. The remaining four points correspond to fibers of type ${\rm I}_1$. 
\item $a^3 = (b - 1)^2 $,  $ b \neq 0$, $ a \neq 0$. In this case, the polynomial ${\rm P}(X)-1$ has 
three distinct roots. However, ${\rm P}(X)+1$ has a root of order two at $-(b-1)/2a$ and a simple 
root at $(b-1)/a$. As in the previous case, the subset $\Sigma$ consists of 5 distinct values. 
$$ \Sigma \ = \ \left \{ \frac{-(b-1)}{2a}, \ \frac{b-1}{a} \right \} \ \cup \ 
\left \{ \ \lambda \ \vert \ {\rm P}(\lambda)=1 \ \right \} $$
The first value in the above list corresponds to a singular fiber of type 
${\rm I}_2$ in $\Theta_2$. The remaining four points correspond to fibers of type ${\rm I}_1$. 
\item $a=0$, $b = \pm 1$. Then 
$$\Sigma \ = \ \{ 0\} \ \cup \ \left \{ \ \frac{1}{\sqrt[3]{2}} \ \theta \ \vert \ \theta^3=b \ \right \}. $$ 
The value $\lambda=0$ corresponds to a singular fiber of type ${\rm I}_3$ in $\Theta_2$. The remaining 
three values of $\Sigma$ correspond to fibers of type ${\rm I}_1$.
\item $a^3=1$, $b = 0$. In this case one has:
$$ {\rm P}(X)-1 \ = \ (2X-a^2)^2(X+a^2), \ \  {\rm P}(X)+1 \ = \ (2X+a^2)^2(X-a^2).$$
Accordingly,
$$\Sigma \ = \ \{ \ \frac{a^2}{2}, \ - \frac{a^2}{2}, -a^2, a^2 \ \}. $$
The first two values in the above list correspond to singular fibers of type ${\rm I}_2$ while the last 
two values correspond to fibers of type ${\rm I}_1$.
\end{itemize}
\end{prop}
\noindent Next, we describe explicitly the involution $\beta$ on ${\rm X}(a,b)$. Note that, in each of the smooth cubics $\Theta_2^{\lambda}$ of $ (\ref{cubics1})$, the point $[1,0,0]$ is an inflection point. If one chooses this point as the origin of the cubic group law on $\Theta_2^{\lambda}$, the point $[0,1,0]$ is a point of order two with 
respect to this law. Moreover, when regarding $\Theta_2^{\lambda}$ as an elliptic fiber in ${\rm X}(a,b)$, one has 
that $[1,0,0]$ and $[0,1,0]$ are the intersections with the two sections $a_1$ and $e_6$. Therefore, the effect of 
$\beta$ on $\Theta_2^{\lambda}$ can be seen, in the coordinates of $(\ref{cubics1})$, as the analytic continuation 
of: 
\begin{equation}
\label{partinv}
\Theta_2^{\lambda} \ \backslash \ \{ [1,0,0], [0,1,0] \} \ \rightarrow \  
\Theta_2^{\lambda} \ \backslash \ \{ [1,0,0], [0,1,0] \} 
\end{equation}
$$ [y,z,w] \ \mapsto \ [-yz,w^2, zw ]. $$
Finally, the full $\beta$ is induced from the analytic involution:
\begin{equation}
\label{inv}
\beta_1 \colon {\rm Q}(a,b) \ \backslash \ \left ( {\rm L}_1 \cup {\rm L}_2 \right ) 
\ \rightarrow \ 
{\rm Q}(a,b) \ \backslash \ 
\left ( {\rm L}_1 \cup {\rm L}_2 \right )
\end{equation}
$$ \beta_1 \left ( [x,y,z,w] \right )   \  \mapsto \ [xz,-yz,w^2,zw] .$$
\subsection{The Elliptic Fibration $\Psi_2$}
\label{step2}
Let ${\rm Y}(a,b)$ be the Kummer surface obtained from ${\rm X}(a,b)$ through the Shioda-Inose 
construction. Recall from section $\ref{proofth1}$ that the alternate fibration $\Theta_2$ 
survives on ${\rm Y}(a,b)$ in the form of a new elliptic fibration $\Psi_2$.   
\par As we already know, the fiber $\Psi_2^{\infty}$ has Kodaira type ${\rm I}^*_{6}$. In this 
section, we describe the location and Kodaira type of the other singular fibers and write an 
explicit formula for the functional invariant ${\rm J}_{\Psi_2}$. 
\par Note that the smooth fibers $\Psi_2^{\lambda} $ are quotients of the cubics 
$\Theta_2^{\lambda}$ of $(\ref{cubics1})$ by the involution $(\ref{partinv})$. By then taking
affine coordinates $[y,z,1]$ on $\Theta_2^{\lambda}$ and defining
$$ u = y^2 - {\rm P}(\lambda) , \ \ v = \frac{1}{2} \ y \left ( z-\frac{1}{z} \right ), $$
one obtains an affine description of $\Psi_2^{\lambda}$ as:
\begin{equation}
v^2 \ =  (u+{\rm P}(\lambda))(u-1)(u+1).
\end{equation}
This can then be easily transformed to a Weierstrass form:
\begin{equation}
v^ 2 \ = \ \left ( u+\frac{1}{3} \ {\rm P}(\lambda) \right )^3 - 
\left ( u+\frac{1}{3} \ {\rm P}(\lambda) \right ) \left ( \frac{1}{3} \ {\rm P}^2(\lambda)  + 1 \right ) + 
\frac{2}{27} \  {\rm P}^3(\lambda) - \frac{2}{3} \  {\rm P}(\lambda) 
\end{equation}
which has as discriminant:
$$ \Delta_{\Psi_2}(\lambda) \ = \ -4 \left ( {\rm P}^2(\lambda)-1 \right )^2.  $$
It follows then that the functional invariant of the elliptic fibration $\Psi_2$ is:
\begin{equation}
\label{jinv22}
\mathcal{J}_{\Psi_2}  (\lambda) \ = \  
\frac{ \left ( 
P^2(\lambda) + 3 
\right )^2}{9 \left  ( P^2(\lambda)-1 \right )^2}.
\end{equation}
As in the previous section, the above information allows us to also describe the location and Kodaira type 
of the singular fibers of $\Psi_2$.
\begin{prop} 
\label{homol1}
The singular fibers of the elliptic fibration $\Psi_2$ on $ {\rm Y}(a,b) $ are located at 
$[1,0]$ (the ${\rm I}_{6}^*$ fiber) and at the points $[\lambda, 1]$ with $\lambda$ belonging to 
the subset:
\begin{equation}
\label{singloc12}
\Sigma : = \ \{ \ \lambda \  \vert \ P(\lambda)^2 = 1 \ \}. 
\end{equation}
The following cases occur:
\begin{itemize}
\item $a^3 \neq (b \pm 1)^2 $. In this case, both polynomials ${\rm P}(X)-1$ and ${\rm P}(X)+1$ have 
three distinct roots. The subset $\Sigma$ consists of six distinct points, each of which corresponds 
to a singular fiber of type ${\rm I}_2$ in $\Theta_2$.
\item $a^3 = (b + 1)^2 $,  $ b \neq 0$, $ a \neq 0$. In this case, ${\rm P}(X)+1$ has three distinct roots. 
However, the polynomial ${\rm P}(X)-1$ has a root of order two at $-(b+1)/2a$ and a simple root at $(b+1)/a$. 
The subset $\Sigma$ consists of 5 distinct points. 
$$ \Sigma \ = \ \left \{ \frac{-(b+1)}{2a}, \ \frac{b+1}{a} \right \} \ \cup \ 
\left \{ \ \lambda \ \vert \ {\rm P}(\lambda)=-1 \ \right \} $$
The first value in the above list corresponds to a singular fiber of type 
${\rm I}_4$ in $\Theta_2$. The remaining four points correspond to fibers of type ${\rm I}_2$. 
\item $a^3 = (b - 1)^2 $,  $ b \neq 0$, $ a \neq 0$. In this case, the polynomial ${\rm P}(X)-1$ has 
three distinct roots. However, ${\rm P}(X)+1$ has a root of order two at $-(b-1)/2a$ and a simple 
root at $(b-1)/a$. As in the previous case, the subset $\Sigma$ consists of 5 distinct values. 
$$ \Sigma \ = \ \left \{ \frac{-(b-1)}{2a}, \ \frac{b-1}{a} \right \} \ \cup \ 
\left \{ \ \lambda \ \vert \ {\rm P}(\lambda)=1 \ \right \} $$
The first value in the above list corresponds to a singular fiber of type 
${\rm I}_4$ in $\Theta_2$. The remaining four points correspond to fibers of type ${\rm I}_2$. 
\item $a=0$, $b = \pm 1$. Then 
$$\Sigma \ = \ \{ 0\} \ \cup \ \left \{ \ \frac{1}{\sqrt[3]{2}} \ \theta \ \vert \ \theta^3=b \ \right \}. $$ 
The value $\lambda=0$ corresponds to a singular fiber of type ${\rm I}_6$ in $\Theta_2$. The remaining 
three values of $\Sigma$ correspond to fibers of type ${\rm I}_2$.
\item $a^3=1$, $b = 0$. In this case one has:
$$ {\rm P}(X)-1 \ = \ (2X-a^2)^2(X+a^2), \ \  {\rm P}(X)+1 \ = \ (2X+a^2)^2(X-a^2).$$
Accordingly,
$$\Sigma \ = \ \{ \ \frac{a^2}{2}, \ - \frac{a^2}{2}, -a^2, a^2 \ \}. $$
The first two values in the above list correspond to singular fibers of type ${\rm I}_4$ while the last 
two values correspond to fibers of type ${\rm I}_2$.
\end{itemize}
\end{prop}

\subsection{A Special Elliptic Fibration on ${\rm Km}({\rm E}_1 \times {\rm E}_2)$}
\label{step3}
\par As we already know from Theorem $\ref{th111}$, the surface ${\rm Y}(a,b)$ can be 
realized in a canonical way as the Kummer surface ${\rm Km}({\rm E}_1 \times {\rm E}_2)$ 
associated to the product of two elliptic curves. Moreover, in this context, 
the elliptic fibration $\Psi_2$ on ${\rm Y}(a,b)$ can be derived directly from the 
Kummer construction. 
 \par Recall from Section $\ref{kummer}$ that the surface ${\rm Km}({\rm E}_1 \times {\rm E}_2)$ has 
a canonical twenty-four curve configuration 
$\{ {\rm H}_i, {\rm G}_j, {\rm E}_{ij} \ \vert \ 0 \leq i,j \leq 3\}$ called the double Kummer pencil. 
\begin{lem}
\label{upsilonf}
Consider the divisor ${\rm D}$ on ${\rm Km}({\rm E}_1 \times {\rm E}_2)$ defined as:
\begin{equation}
\label{specialdiv}
{\rm D} \ = \ {\rm E}_{21} + {\rm E}_{31} + 2 \left ( {\rm G}_1 + {\rm E}_{01} + {\rm H}_0 + 
{\rm E}_{00} + {\rm G}_0 + {\rm E}_{10} + {\rm H}_1 \right ) + {\rm E}_{12} + {\rm E}_{13}. 
\end{equation}
Then ${\rm D}^2=0$ and $ \vert D \vert $ is a pencil inducing an elliptic fibration 
$ \Upsilon_2 \colon {\rm Km}({\rm E}_1 \times {\rm E}_2) \rightarrow \mathbb{P}^1 $. The 
divisor ${\rm D}$ is a singular fiber for $\Upsilon_2$ and has Kodaira type ${\rm I}_6^*$. 
The four smooth rational curves ${\rm H}_2,{\rm H}_3, {\rm G}_2$ and ${\rm G}_3$ form four disjoint sections of 
$\Upsilon_2$. 
\end{lem}
 $$ 
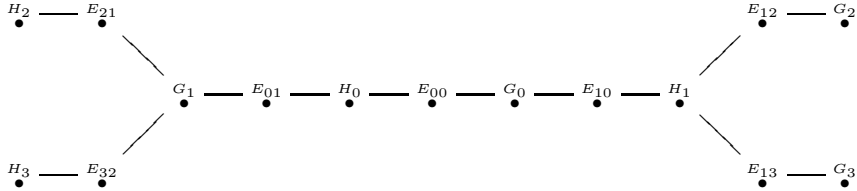
\def\objectstyle{\scriptstyle}
\def\labelstyle{\scriptstyle}
\xymatrix @-0.8pc  {
\stackrel{H_2}{\bullet} \ar @{-} [r] & \stackrel{E_{21}}{\bullet} \ar @{-} [dr] & & & & &  & & & \stackrel{E_{12}}{\bullet} \ar @{-} [dl] 
 & \stackrel{G_{2}}{\bullet} \ar @{-} [l] \\
& & \stackrel{G_1}{\bullet} \ar @{-} [r] \ar @{-} [dl] &
\stackrel{E_{01}}{\bullet} \ar @{-} [r] &
\stackrel{H_0}{\bullet} \ar @{-} [r] &
\stackrel{E_{00}}{\bullet} \ar @{-} [r] &
\stackrel{G_0}{\bullet} \ar @{-} [r] &
\stackrel{E_{10}}{\bullet} \ar @{-} [r] &
\stackrel{H_{1}}{\bullet} \ar @{-} [dr] & \\
\stackrel{H_3}{\bullet} \ar @{-} [r] & \stackrel{E_{32}}{\bullet} & & & & &   & & & \stackrel{E_{13}}{\bullet} & \stackrel{G_{3}}{\bullet} \ar @{-} [l]\\
}
$$
\bpf
The above assertion is a consequence of a classical theorem due to 
I.I. Pjatecki\u{i}-\v{S}apiro and I.R. \v{S}afarevi\v{c} (\cite{shapiro}, Chapter 3, Theorem 1). 
\epf
\begin{rem}
A different selection of the double Kummer pencil curves defining the divisor $(\ref{specialdiv})$ alters 
the elliptic fibration $\Upsilon_2$ by an analytic automorphism of ${\rm Km}({\rm E}_1 \times {\rm E}_2)$. 
The equivalence class of $\Upsilon_2$ is therefore well-defined. 
\end{rem} 
\begin{rem}
In \cite{oguiso}, K. Oguiso classified all jacobian fibrations on a Kummer surface associated to 
a product of two non-isogenous elliptic curves. The elliptic fibration $\Upsilon_2$ defined above 
appears as $\mathcal{J}_5$ in Oguiso's classification. It is the only jacobian fibration on such a 
surface that admits a singular fiber of Kodaira type ${\rm I}^*_{6}$.  
\end{rem}
\noindent By virtue of the geometric correspondence 
$$ {\rm X}(a,b) \ \longrightarrow \ {\rm E}_1 \times {\rm E}_2 $$
described in the first part of the paper, one has, as an intermediate step, an isomorphism 
$${\rm Y}(a,b) \ \simeq \ {\rm Km}({\rm E}_1 \times {\rm E}_2)$$ 
that maps the jacobian fibration 
$\Psi_2$ on ${\rm Y}(a,b)$ to the jacobian fibration $\Upsilon_2$ on ${\rm Km}({\rm E}_1 \times {\rm E}_2)$.
\par This fact allows one to realize an explicit relation between the Inose parameters $a,b$ of 
the ${\rm M}$-polarized K3 surface ${\rm X}(a,b)$ and the ${\rm J}$-invariants of the two resulting 
elliptic curves ${\rm E}_1$ and ${\rm E}_2$. 
\noindent In light of Corollary $\ref{cormain}$, the two jacobian fibrations $\Psi_2$ on ${\rm Y}(a,b)$ 
and $\Upsilon_2$ on ${\rm Km}({\rm E}_2 \times {\rm E}_2)$ are equivalent if and only if 
their functional invariant and singular locus differ by a projective transformation and the Kodaira 
types of their singular fibers match. 
\par We have already described in detail the functional invariant and the location and type of 
the singular fibers of $\Psi_2$. In what follows we shall perform a similar analysis for $\Upsilon_2$. 
The comparison between these two pieces of data will then allow us to prove the main statement of 
Theorem $\ref{main}$.
\begin{cl}
\label{lastclaim}
The two elliptic fibrations $\Psi_2$ and $\Upsilon_2$ have equivalent functional and homological 
invariants if and only if 
$$ {\rm J}({\rm E}_1) + {\rm J}({\rm E}_2) \ = \ a^3-b^2+1, \ \ \ \ 
{\rm J}({\rm E}_1) \cdot {\rm J}({\rm E}_2) \ = \ a^3.$$ 
\end{cl}
\subsection{Description of the Elliptic Fibration $\Upsilon_2$ on ${\rm Km}({\rm E}_1 \times {\rm E}_2)$}
\label{adouafib}  
\par It is a standard fact that any given elliptic curve can be realized as a projective Legendre cubic
$$ \{ \ y^2w = x(x-w)(x-\lambda w) \ \} \ \subset \ \mathbb{P}^2 $$
for some $\lambda \in \Cee \setminus \{ 0,1\}$. We shall assume therefore that $\alpha, \beta \in \Cee $ are chosen 
such that ${\rm E}_1$ and ${\rm E}_2$ are isomorphic with the above cubics for $\lambda=\alpha$ and 
$\lambda = \beta$, respectively. The ${\rm J}$-invariants of the two curves can then be computed as:
$$ {\rm J}({\rm E}_1) \ = \ \frac{4(\alpha^2-\alpha+1)^3}{27\alpha^2(\alpha-1)^2}, \ \ \ 
{\rm J}({\rm E}_2) \ = \ \frac{4(\beta^2-\beta+1)^3}{27\beta^2(\beta-1)^2}. $$ 
In this context, an explicit model for the Kummer surface ${\rm Km}({\rm E}_1 \times {\rm E}_2) $ 
can be constructed (see \cite{cassels, inose}) by taking the minimal resolution of the quartic surface:
\begin{equation}
\label{kummerquartic}
\{ \ z^2xy \ = (x-w)(x-\alpha w)(y-w)(y - \beta w) \ \} \ \subset \ \mathbb{P}^3.
\end{equation} 
Note that, generically, the quartic surface $(\ref{kummerquartic})$ has seven rational double point singularities located 
at:
$$ 
[1,0,0,0], \ \ 
[0,0,1,0], \ \ 
[0,0,1,0]
$$
$$ 
[1,1,0,1], \ \ 
[\alpha, 1,0,1], \ \ 
[1,\beta,0,1], \ \ 
[\alpha, \beta, 0, 1].
$$
The first three are rational double points of type ${\rm A}_3$. The last four are singularities 
of type ${\rm A}_1$. One can therefore reconstruct the double Kummer pencil on the minimal 
resolution of $(\ref{kummerquartic})$ by taking:
$$ {\rm H}_0 + {\rm E}_{00} + {\rm G}_0 \ \ = \ \ {\rm A}_3 \ {\rm configuration \ associated \ to} \ [0,0,1,0]$$
$$ {\rm E}_{21}+{\rm G}_1+{\rm E}_{31} \ \ = \ \ {\rm A}_3 \ {\rm configuration \ associated \ to} \ [1,0,0,0]$$
$$ {\rm E}_{12}+{\rm H}_1+{\rm E}_{13} \ \ = \ \ {\rm A}_3 \ {\rm configuration \ associated \ to} \ [0,1,0,0]$$
$$ {\rm E}_{22} \ \ = \ \ {\rm A}_1 \ {\rm curve \ associated \ to} \ [1,1,0,1] $$
$$ {\rm E}_{32} \ \ = \ \ {\rm A}_1 \ {\rm curve \ associated \ to} \ [\alpha,1,0,1] $$
$$ {\rm E}_{23} \ \ = \ \ {\rm A}_1 \ {\rm curve \ associated \ to} \ [1,\beta,0,1] $$
$$ {\rm E}_{33} \ \ = \ \ {\rm A}_1 \ {\rm curve \ associated \ to} \ [\alpha,\beta,0,1] $$
$$ {\rm H}_2 \ \ = \ \ {\rm proper \ transform \ of} \ \{ x=w, \ z=0 \}$$
$$ {\rm H}_3 \ \ = \ \ {\rm proper \ transform \ of} \ \{ x=\alpha w, \ z=0 \}$$
$$ {\rm G}_2 \ \ = \ \ {\rm proper \ transform \ of} \ \{ y=w, \ z=0 \}$$
$$ {\rm G}_3 \ \ = \ \ {\rm proper \ transform \ of} \ \{ y=\beta w, \ z=0 \}$$
$$ {\rm E}_{01} \ \ = \ \ {\rm proper \ transform \ of} \ \{ x=w=0 \}$$
$$ {\rm E}_{10} \ \ = \ \ {\rm proper \ transform \ of} \ \{ y=w=0 \}$$
$$ {\rm E}_{11} \ \ = \ \ {\rm proper \ transform \ of} \ \{ w=0, \ z^2=xy \}$$
$$ {\rm E}_{02} \ \ = \ \ {\rm proper \ transform \ of} \ \{ x=0, \ y=w \}$$
$$ {\rm E}_{03} \ \ = \ \ {\rm proper \ transform \ of} \ \{ x=0, \ y=\beta w \}$$
$$ {\rm E}_{20} \ \ = \ \ {\rm proper \ transform \ of} \ \{ y=0, \ x=w \}$$
$$ {\rm E}_{30} \ \ = \ \ {\rm proper \ transform \ of} \ \{ y=0, \ x=\alpha w \}$$
A simple analysis of the curves ${\rm H}_i$, ${\rm G}_j$ and locations of the intersections 
with ${\rm E}_{ij}$ allows one to conclude that the minimal resolution of 
$(\ref{kummerquartic})$ is canonically 
isomorphic to ${\rm Km}({\rm E}_1 \times {\rm E}_2)$. 
\begin{rem}
The birational morphism ${\rm Km}({\rm E}_1 \times {\rm E}_2) \rightarrow \mathbb{P}^3$ whose 
image is the quartic surface $(\ref{kummerquartic})$ can also be defined directly from the 
double Kummer pencil by taking the projective morphism associated to the base-point free 
linear system $\vert {\rm V} \vert $ given by:
$$ {\rm V} \ = \ {\rm E}_{12} + 2 {\rm H}_1 + {\rm E}_{13} + {\rm E}_{10} + {\rm G}_{0} + 
{\rm E}_{00} + 
{\rm H}_0 + {\rm E}_{01} + 2{\rm G}_1 + {\rm E}_{21} + {\rm E}_{31} + 2 {\rm E}_{11}. $$  
$$ \def\objectstyle{\scriptstyle}
\def\labelstyle{\scriptstyle}
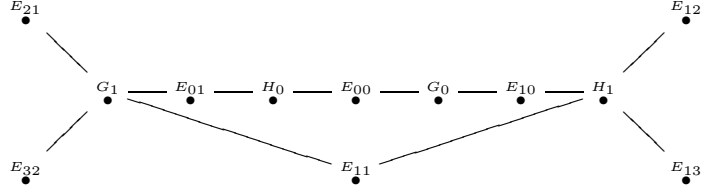
\xymatrix @-0.8pc  {
\stackrel{E_{21}}{\bullet} \ar @{-} [dr] & & & & &  & & & \stackrel{E_{12}}{\bullet} \ar @{-} [dl] 
  \\
& \stackrel{G_1}{\bullet} \ar @{-} [r] \ar @{-} [drrr] \ar @{-} [dl] &
\stackrel{E_{01}}{\bullet} \ar @{-} [r] &
\stackrel{H_0}{\bullet} \ar @{-} [r] &
\stackrel{E_{00}}{\bullet} \ar @{-} [r] &
\stackrel{G_0}{\bullet} \ar @{-} [r] &
\stackrel{E_{10}}{\bullet} \ar @{-} [r] &
\stackrel{H_{1}}{\bullet} \ar @{-} [dr]  \ar @{-} [dlll] \\
 \stackrel{E_{32}}{\bullet} & & & & \stackrel{E_{11}}{\bullet}&   & & & \stackrel{E_{13}}{\bullet} \\
}
$$
\end{rem}
\noindent The advantage of the realization of ${\rm Km}({\rm E}_1 \times {\rm E}_2)$ as the 
quartic surface $(\ref{kummerquartic})$, follows from the fact that, in this context, one can 
explicitly construct the jacobian fibration $\Upsilon_2$ of Lemma $\ref{upsilonf}$. 
This elliptic fibration is induced by the rational map:
$$ [x,y,z,w] \ \mapsto \ \left [ {\rm R}(x,y,w), \ xy \right ]$$
where ${\rm R}(x,y,w)$ is the quadratic polynomial:
$$ {\rm R}(x,y,w) \ = \ \left ( - \frac{1}{\alpha } \right ) x^2 +  \left ( - \frac{1}{\beta } \right ) y^2 + 
\left ( \frac{\alpha + 1}{\alpha } \right ) xw +  \left (  \frac{\beta + 1}{\beta } \right )yw - w^2. $$
The ${\rm I}^*_6$ fiber of $\Upsilon_2$ appears over the point $[1,0]$. Away from this location, the generic 
smooth elliptic fiber 
$$\Upsilon_2^{\lambda} : = \Upsilon_2^{-1}\left ( [\mu,1] \right ) $$ 
can be regarded as the double cover of the projective conic in $\mathbb{P}^2(x,y,w)$:
\begin{equation}
\label{conic}
{\rm R}(x,y,w) \ = \ \mu xy 
\end{equation}
branched at the four points:
\begin{equation}
\label{branchpoints}
 [1,(1-\mu)\beta + 1, 1 ], \ \ [\alpha,(1-\mu \alpha )\beta + 1, 1 ], 
 \end{equation}
$$ [(1-\mu)\alpha + 1, 1, 1 ], \ \ [(1-\mu \beta )\alpha + 1, \beta, 1 ]. $$
One encounters singular fibers if the conic $(\ref{conic})$ is singular, or if at least two 
of the above four branch points coincide.  This argument allows one to conclude that the 
points on the base of the fibration $\Upsilon_2$ associated to singular fibers 
(away from the ${\rm I}^*_6$ fiber) are of type $ [\mu,1]$ with 
$\mu$ belonging to the set:
$$ \Sigma_{\Upsilon_2} : =  \left \{ \ 1, \ \frac{1}{\alpha}, \ \frac{1}{\beta }, \ \frac{1}{\alpha \beta }, \ 
\frac{\alpha \beta + 1}{\alpha \beta }, \ \frac{\alpha + \beta}{\alpha \beta } \ \right \}.$$
\begin{rem}
For generic choices of $\alpha$ and $\beta$, the above set contains six distinct points, each of which 
determines an ${\rm I}_2$ fiber in $\Upsilon_2$. However, it may happen that two, or more, of the above 
six values coincide. This happens precisely when:
\begin{equation}
\label{singg}
\beta \ \in \ \left \{ \alpha, \ \frac{1}{\alpha}, \ 1-\alpha, \ \frac{1}{1-\alpha}, \ \frac{\alpha}{\alpha-1}, 
\ \frac{\alpha-1}{\alpha} \  \right \}, 
\end{equation}
a condition that is also equivalent to ${\rm J}({\rm E}_1)={\rm J}({\rm E}_2)$.
\end{rem}
\begin{lem}
The functional invariant of the jacobian fibration $\Upsilon_2$ has the form:
\begin{equation}
\label{funcinvupsilon}
\mathcal{J}_{\Upsilon_2}(\mu) \ = \ 
\frac{
4 \left ( \alpha^4\beta^4 \ {\rm D}(\mu) + (\alpha-1)^2(\beta-1)^2 \right )^3
}
{
27 \ \alpha^8\beta^8 (\alpha-1)^4(\beta-1)^4 \ {\rm D}^2(\mu)
}
\end{equation}
where:
$$ {\rm D}(\mu) := \ \left ( \mu - 1 \right )
\left ( \mu - \frac{1}{\alpha  } \right )
\left ( \mu - \frac{1}{\beta } \right )
\left ( \mu - \frac{1}{\alpha \beta } \right )
\left ( \mu - \frac{\alpha \beta+ 1}{\alpha \beta } \right )
\left ( \mu - \frac{\alpha + \beta}{\alpha \beta } \right ).$$
\end{lem}
\bpf
We accomplish the computation of $\mathcal{J}_{\Upsilon_2}(\mu)$ through the 
following sequence of steps.  
\begin{enumerate}
\item Construct an explicit isomorphism $i_{\mu} $ between the conic $(\ref{conic})$ and $\mathbb{P}^1$.
\item Perform a projective automorphism of $\mathbb{P}^1$ such that the images through $i_{\mu}$ of 
the four branch points $(\ref{branchpoints})$ are sent to $[0,1],[1,1],[r,1]$ and $[1,0]$.  
\item Evaluate $\mathcal{J}_{\Upsilon_2}(\mu)$ as:
\begin{equation}
\frac{
4(r^2-r+1)^3
}
{
27r^2(r-1)^2
}.
\end{equation}
\end{enumerate}
In order to complete the first step, let us note that:
$$ {\rm R}(x,y,w) - \mu xy \ = \ $$
$$ - \left [ w - \left ( \frac{\alpha + 1}{2 \alpha } \right ) x - \left ( \frac{\beta + 1}{2 \beta } \right ) y \right ]^2 + \left ( \frac{\alpha - 1}{2 \alpha } \right )^2 x^2 + \left ( \frac{\beta - 1}{2 \beta } \right )^2 y^2 
+ \left ( \frac{(\alpha+1)(\beta+1)}{2 \alpha \beta } - \mu \right ) xy \ =  $$
$$ - \left [ 
w - \left ( \frac{\alpha + 1}{2 \alpha } \right ) x - \left ( \frac{\beta + 1}{2 \beta } \right ) y 
\right ]^2 +
\left [ 
\left ( \frac{\alpha - 1}{2 \alpha } \right ) x + 
\left ( \frac{(\alpha+1)(\beta+1)}{2 \alpha \beta } - \mu \right ) \left ( \frac{ \alpha}{\alpha-1} \right )y
\right ]^2 \ + 
$$
$$
+ \ 
\left [ 
\left ( \frac{\beta - 1}{2 \beta } \right )^2 - \left ( \frac{(\alpha+1)(\beta+1)}{2 \alpha \beta } - \mu \right )^2 \left ( \frac{ \alpha}{\alpha-1} \right )^2
\right ] y^2 \ = 
$$
$$ - \left [ 
w - \left ( \frac{\alpha + 1}{2 \alpha } \right ) x - \left ( \frac{\beta + 1}{2 \beta } \right ) y 
\right ]^2 +
\left [ 
\left ( \frac{\alpha - 1}{2 \alpha } \right ) x + 
\left ( \frac{(\alpha+1)(\beta+1)}{2 \alpha \beta } - \mu \right ) \left ( \frac{ \alpha}{\alpha-1} \right )y
\right ]^2 \ + 
$$
$$
+ \ 
\left [ 
\left ( \frac{(\alpha-1)(\beta - 1)}{2 \alpha \beta } \right )^2 - \left ( \frac{(\alpha+1)(\beta+1)}{2 \alpha \beta } - \mu \right )^2 
\right ] \left ( \frac{ \alpha}{\alpha-1} \right )^2 y^2 \ = 
$$
$$ - \left [ 
w - \left ( \frac{\alpha + 1}{2 \alpha } \right ) x - \left ( \frac{\beta + 1}{2 \beta } \right ) y 
\right ]^2 +
\left [ 
\left ( \frac{\alpha - 1}{2 \alpha } \right ) x + 
\left ( \frac{(\alpha+1)(\beta+1)}{2 \alpha \beta } - \mu \right ) \left ( \frac{ \alpha}{\alpha-1} \right )y
\right ]^2 \ + 
$$
$$
- \ 
\left [ 
\left ( \mu - \frac{\alpha + \beta }{\alpha \beta } \right ) 
\left ( \mu - \frac{\alpha \beta + 1 }{\alpha \beta }\right ) 
\right ] \left ( \frac{ \alpha}{\alpha-1} \right )^2 y^2 \ \ = \ 
$$
$$ = \ 
\left ( w - \frac{1}{\alpha} \ x - \left ( \mu - \frac{\beta + 1}{\alpha \beta} \right ) \left ( \frac{\alpha}{\alpha-1}\right )y\right ) \left ( -w + x - \left ( \mu - \frac{\beta + 1}{\beta} \right ) \left ( \frac{\alpha}{\alpha-1}\right )y \right ) \ - \ 
$$
$$
- \ 
\left [ 
\left ( \mu - \frac{\alpha + \beta }{\alpha \beta } \right ) 
\left ( \mu - \frac{\alpha \beta + 1 }{\alpha \beta }\right ) 
\right ] \left ( \frac{ \alpha}{\alpha-1} \right )^2 y^2 \
$$
The change in projective coordinates:
$$ x_1 \ = \ -w + x - \left ( \mu - \frac{\beta + 1}{\beta} \right ) \left ( \frac{\alpha}{\alpha-1}\right )y ,
 \ \ \ \ \ y_1 \ = \ \left ( \frac{\alpha}{\alpha-1}\right )y $$
$$ w_1 \ = \  w - \frac{1}{\alpha} \ x - \left ( \mu - \frac{\beta + 1}{\alpha \beta} \right ) \left ( \frac{\alpha}{\alpha-1}\right )y $$
allows one to rewrite the conic $ (\ref{conic})$ as:
$$ x_1 w_1 = \Delta  y_1^2 $$
where:
$$ \Delta \ = \  \left ( \mu - \frac{\alpha + \beta }{\alpha \beta } \right ) 
\left ( \mu - \frac{\alpha \beta + 1 }{\alpha \beta }\right ). 
$$
This yields the parametrization of $(\ref{conic})$ via the embedding:
$$ \mathbb{P}^1 \hookrightarrow \mathbb{P}^2, \ \ [u,v] \mapsto [ u^2,uv,\Delta v^2] $$
with the inverse map $i_{\mu}$ given by the analytic continuation of $[x_1,y_1,w_1] \mapsto [x_1,y_1]$. 
This procedure results in an identification between the conic $(\ref{conic})$ and $ \mathbb{P}^1$ which 
sends the four branch points $(\ref{branchpoints})$ to
$$ [\beta + 1-\beta \mu, \ \beta], \ \ \ 
[\alpha \beta \Delta , \ \beta+ 1 - \mu \alpha \beta], $$
$$ [\alpha \beta + 1 - \mu \alpha \beta , \ \beta ], \ \ \ 
[\alpha + \beta - \mu \alpha \beta , \beta ].$$
In accordance with the plan presented earlier, we take:
$$ r \ = \ \frac{
\left ( \left ( \alpha \beta + 1 - \mu \alpha \beta \right ) - \left ( \beta + 1-\beta \mu \right ) \right ) 
\left ( \left ( \alpha \beta \Delta\right ) - \left ( \alpha + \beta - \mu \alpha \beta\right ) \right ) 
}{
\left ( \left ( \alpha \beta + 1 - \mu \alpha \beta\right ) - \left ( \alpha + \beta - \mu \alpha \beta\right ) \right )
\left ( \left ( \alpha \beta \Delta \right ) - \left ( \beta + 1-\beta \mu \right ) \right )
} \ = \ $$
$$ 
= \ \frac{(\mu-1)(\mu \alpha \beta-1)(\mu \alpha \beta - \alpha -\beta)}{
(\alpha-1)(\beta-1)}.
$$
Then the functional invariant is computed as:
$$
\mathcal{J}_{\Upsilon_2}(\mu) \ = \ \frac{4(r^2-r+1)^3}{27r^2(r-1)^2} \ = \ 
\frac{
4 \left ( \alpha^4\beta^4 \ {\rm D}(\mu) + (\alpha-1)^2(\beta-1)^2 \right )^3
}
{
27 \ \alpha^8\beta^8 (\alpha-1)^4(\beta-1)^4 \ {\rm D}^2(\mu)
}.
$$
\epf \

\noindent The above discussion also provides the homological invariant data of the fibration 
$\Upsilon_2$.
\begin{cor}
In addition to the ${\rm I}_6^*$ singular fiber which appears over the point $[1,0]$, the elliptic 
fibrations $\Upsilon_2$ has singular fibers at the points $ [\mu, 1 ]$ with $ \mu $ belonging 
to the set:
$$ \Sigma_{\Upsilon_2} \ = \ \{ \mu \ \vert \ {\rm D}(\mu)=0 \ \}. $$  
The following cases can occur:
\begin{itemize}
\item [(a)] ${\rm J}({\rm E}_1) \neq {\rm J}({\rm E}_2).$ In this case $\Sigma_{\Upsilon_2}$ has six 
distinct points and each of them corresponds to an ${\rm I}_2$ singular fiber.
\item [(b)] ${\rm J}({\rm E}_1) = {\rm J}({\rm E}_2) \notin \{ 0,1\} $. In this case the polynomial 
${\rm D}(\mu)$ has five distinct roots, one of which is of order two. The order-two root corresponds 
to a singular fiber of type ${\rm I}_4$. The remaining four roots correspond to ${\rm I}_2$ fibers. 
\item [(c)] ${\rm J}({\rm E}_1) = {\rm J}({\rm E}_2) = 1. $ In this case the polynomial 
${\rm D}(\mu)$ has four distinct roots, two of which have order two. The two roots of order two 
correspond to singular fibers of type ${\rm I}_4$. The remaining two roots correspond to fibers 
of type ${\rm I}_2$. 
\item [(d)] ${\rm J}({\rm E}_1) = {\rm J}({\rm E}_2) = 0. $ In this case the polynomial 
${\rm D}(\mu)$ has four distinct roots, one of which is of order three. The order-three root 
corresponds to singular fiber of type ${\rm I}_3$. The remaining three roots correspond to fibers 
of type ${\rm I}_2$. 
\end{itemize}
\end{cor}
\subsection{Proof of Claim $\ref{lastclaim}$}
\noindent Recall the analysis of Sections $\ref{step2}$ and $\ref{adouafib}$. Both fibrations $\Psi_2$ and $\Upsilon_2$ have the ${\rm I}_6^*$ singular fiber 
located over the point $ [1,0]$ and their respective functional invariants, as described in $(\ref{jinv22})$ 
and $(\ref{funcinvupsilon})$, are :
\begin{equation}
\label{jinvpsi}
\mathcal{J}_{\Psi_2}  (\lambda) \ = \  
\frac{ \left ( 
{\rm P}^2(\lambda) + 3 
\right )^2}{9 \left  ( {\rm P}^2(\lambda)-1 \right )^2}, \ \ \ \ {\rm P}(\lambda) = 4\lambda^3 -3a \lambda - b.
\end{equation}   
  \begin{equation}
\label{jinvupsilon}
\mathcal{J}_{\Upsilon_2}(\mu) \ = \ 
\frac{
4 \left ( \alpha^4\beta^4 \ {\rm D}(\mu) + (\alpha-1)^2(\beta-1)^2 \right )^3
}
{
27 \ \alpha^8\beta^8 (\alpha-1)^4(\beta-1)^4 \ {\rm D}^2(\mu)
}.
\end{equation}
The main polynomial in the denominator of ${\rm J}_{\Psi_2}(\lambda)$ is ${\rm P}^2(\lambda)-1$. Its (generic) six roots are naturally divided into two sets of three roots, each three-set having the sum of its elements equal to zero. A similar feature can 
be observed in the denominator of ${\rm J}_{\Upsilon_2}(\mu)$. The main polynomial present there is 
${\rm D}(\mu)$ whose (generic) six roots can be partitioned into two sets of three with identical sum.
\begin{equation}
\label{tt123}
\left \{ 
\ 1, \ \frac{1}{\alpha}, \ \frac{1}{\beta }, \ \frac{1}{\alpha \beta }, \ 
\frac{\alpha \beta + 1}{\alpha \beta }, \ \frac{\alpha + \beta}{\alpha \beta } \
\right \} 
\ = 
 \ 
\left \{ 
\ 1, \ \frac{1}{\alpha \beta }, \ \frac{\alpha + \beta}{\alpha \beta } \
\right \} \ \ \cup \ \  
 \left \{ 
\ \frac{1}{\alpha}, \ \frac{1}{\beta }, \ 
\frac{\alpha \beta + 1}{\alpha \beta } \
\right \} 
\end{equation}
The two fibrations $\Psi_2$ and $\Upsilon_2$ have equivalent functional invariant and homological invariant 
data if and only if there exists an invertible affine transformation $\Xi(\lambda) = q \lambda + p$ with 
$p,q\in \Cee$ ($q \neq 0$) such that 
\begin{equation}
\label{eer}
\mathcal{J}_{\Psi_2}(\lambda) = \mathcal{J}_{\rm \Upsilon_2} (\Xi(\lambda)) 
\end{equation}
and $\Xi$ sends the roots of ${\rm P}(\lambda) \pm 1$ to the two subsets in $(\ref{tt123})$ while 
preserving the 
homological type.   
\par As a first observation, we note that it follows that $p = (\alpha+1)(\beta+1)/3\alpha \beta $. Let then:
$$ {\rm D}_1(\mu) \ = \ (\mu - 1) \left ( \mu - \frac{1}{\alpha \beta} \right ) 
\left (\mu - \frac{\alpha+ \beta}{\alpha \beta} \right )$$
$$ {\rm D}_2(\mu) \ = \ \left ( \mu - \frac{1}{\alpha} \right ) \left ( \mu - \frac{1}{\beta } \right ) 
\left (\mu - \frac{\alpha\beta+1}{\alpha \beta} \right ).$$
Two possibilities can occur:
\begin{itemize}
\item [(a)] $q^3 \left ( {\rm P}(\lambda)-1 \right )  \ = \ 4{\rm D}_1(q\lambda +p)$ and 
$q^3 \left ( {\rm P}(\lambda)+1 \right )  \ = \ 4{\rm D}_2(q\lambda +p)$
\item [(b)] $q^3 \left ( {\rm P}(\lambda)-1 \right )  \ = \ 4{\rm D}_2(q\lambda +p)$ and 
$q^3 \left ( {\rm P}(\lambda)+1 \right )  \ = \ 4{\rm D}_1(q\lambda +p).$
\end{itemize}
Making the constant terms coincide implies, in the two cases:
$$ b \ = \ \pm \ \frac{(\alpha-2)(\alpha +1)(2\alpha-1)(\beta-2)(\beta +1)(2\beta-1)}
{27 \alpha(\alpha-1)\beta(\beta-1)} \ $$
and 
\begin{equation}
\label{ccaa}
 q^3 \ = \ - \frac{2(\alpha-1)(\beta-1)}{\alpha^2\beta^2}.
 \end{equation}
The first equality, in turn, requires the a priori condition.
\begin{equation}
\label{cond111}
b^2 \ = \ \ \frac{(\alpha-2)^2(\alpha +1)^2(2\alpha-1)^2(\beta-2)^2(\beta +1)^2(2\beta-1)^2}
{729 \alpha^2(\alpha-1)^2\beta^2(\beta-1)^2}\ = \ 
\left ( {\rm J}({\rm E}_1) -1 \right ) \left ( {\rm J}({\rm E}_2) -1 \right ).
\end{equation}
Continuing the argument, we note that when imposing the equality of the linear terms in the above cases, one is 
led to:
\begin{equation}
 a \ = \ \frac{4(\alpha^2-\alpha+1)(\beta^2-\beta+1)}{9 \alpha^2 \beta^2 q^2}.
 \end{equation}
This constraint, in connection with $(\ref{ccaa})$, requires a second a priori condition:
\begin{equation}
\label{cond222}
a^3 \ = \ \frac{16(\alpha^2-\alpha+1)^3(\beta^2-\beta+1)^3}{729\alpha^2(\alpha-1)^2\beta^2(\beta-1)^2} \ 
= \ {\rm J}({\rm E}_1) \cdot {\rm J}({\rm E}_2).  
\end{equation}
Then, depending on the case in question, one can eventually solve for $q$, yielding:
$$ q \ = \ \pm \frac{-9(\alpha-1)(\beta-1)}{2(\alpha^2-\alpha+1)(\beta^2-\beta+1)}.$$
The equivalence $(\ref{eer})$ of functional invariants is immediately verified.
To finish the proof of Claim $\ref{lastclaim}$, we note that conditions $(\ref{cond111})$ and 
$(\ref{cond222})$ are equivalent to:
$$ {\rm J}({\rm E}_1) + {\rm J}({\rm E}_2)\ = \ a^3-b^2+1, \ \ \ 
{\rm J}({\rm E}_1) \cdot {\rm J}({\rm E}_2) \ = \ a^3. $$  

\section{A String Duality Point of View}
\label{finalcomments}
Following the works 
of Vafa \cite{vafa96} and Sen \cite{sen96} in 1996, 
it was noted that the geometry underlying elliptic $K3$ surfaces with 
section is related to the geometry of elliptic curves endowed with certain 
flat principal $G$-bundles and an additional parameter called the B-field. 
This non-trivial connection appears in string theory 
as the eight-dimensional 
manifestation of the phenomenon called F-theory/heterotic string duality. Over the past ten years 
the correspondence has been analyzed extensively (\cite{curio, clingher1}) from a purely 
mathematical point of view. As it turns out, it leads to a beautiful 
geometric picture which links together moduli spaces for 
these two seemingly distinct types of geometrical objects: elliptic 
$K3$ surfaces with section and flat bundles over elliptic curves. 
\par In brief, what happens 
is the following. On the F-theory side, one has the moduli space $\mathcal{M}_{{\rm K3}}$ of elliptic 
K3 surfaces with section. This is a quasi-projective analytic variety of complex dimension eighteen. 
However, $\mathcal{M}_{{\rm K3}}$ is not compact. Nevertheless, there exists a nice smooth partial 
compactification $\mathcal{M}_{{\rm K3}} \subset \overline{\mathcal{M}}_{{\rm K3}}$, consisting of 
an enlargement of the original space by adding two Type II Mumford boundary divisors $\mathcal{D}_1$ and $\mathcal{D}_2$. Geometrically, the points of the two compactifying divisors correspond to 
Type II stable elliptic K3 surfaces. These are special degenerations of K3 surfaces realized as a union 
${\rm V}_1 \cup {\rm V}_2$ of two rational surfaces meeting over a common elliptic curve ${\rm E}$ 
which is anti-canonical on both ${\rm V}_1$ and ${\rm V}_2$. 
\par On the heterotic side, one has to consider two moduli spaces $\mathcal{M}_{{\rm het}}^{\rm G}$ of triples 
$({\rm E}, {\rm P}, {\rm B})$ consisting of elliptic curves, flat ${\rm G}$-bundles, and B-fields. There are two choices of Lie groups ${\rm G}$:
\begin{equation}
\label{picture}
\eeight \ \ \ \ \ \ \spint.
\end{equation}
The moduli space $\mathcal{M}_{{\rm E}, {\rm G}}$ associated to the first two components 
$({\rm E}, {\rm P})$ of the above triples is (as described in \cite{friedmanmorgan}) 
a quasi-projective analytic space of 
complex dimension seventeen. The actual heterotic moduli space 
$ \mathcal{M}_{{\rm het}}^{{\rm G}}$ (described in \cite{clingher2}) fibers naturally as a 
holomorphic $\mathbb{C}^*$ fibration over 
$\mathcal{M}_{{\rm E}, {\rm G}}$.
\par In this context, the mathematical facts underlying the string duality can be summarized as follows. 
Each of the two moduli spaces $\mathcal{M}_{{\rm E}, {\rm G}}$ 
associated to the two choices of possible Lie groups on the heterotic side 
is naturally isomorphic to one of the corresponding Type II Mumford boundary divisors $\mathcal{D}_1$ and 
$\mathcal{D}_2$ from the F-theory side. Moreover, there exists a holomorphic identification between an open subset 
of $\mathcal{M}_{{\rm het}}^{{\rm G}}$ neighboring the cusps of the $\mathbb{C}^*$-fibration 
over $\mathcal{M}_{{\rm E},{\rm G}}$ and a special subset (of large complex structures) 
of $\mathcal{M}_{{\rm K3}}$ which makes an open neighborhood of the corresponding boundary 
divisor. We refer the reader to \cite{clingher3, clingher1} for further details and proofs.  
\par The above holomorphic identification between the appropriate regions of 
$\mathcal{M}_{{\rm K3}}$ and ${\rm M}_{{\rm het}}^{{\rm G}}$ is however defined Hodge-theoretically and 
therefore is not fully satisfactory from a geometer's point of view. As with any string duality, 
one would like to have a purely geometrical pattern that connects the spaces and structures
appearing on the two sides of the duality correspondence. 
\par Such a geometric connection, in the context 
of the above duality, has been known for some time, but only in the stable limit, i.e. on the boundary 
of the moduli spaces \cite{aspinwall, fmw}. As mentioned earlier, on the F-theory side 
this limit corresponds to stable K3 surfaces, whereas on the heterotic side it corresponds to 
${\rm B}=0$. Given a Type II stable K3 surface ${\rm V}_1 \cup {\rm V}_2$, 
one can obtain the heterotic elliptic curve ${\rm E}$ by just taking the common curve 
${\rm V}_1 \cap {\rm V}_2$ and the heterotic ${\rm G}$-bundle can also be derived explicitly 
from the geometry of the rational surfaces ${\rm V}_1$ and ${\rm V}_2$.   
\par It is natural to ask whether there exists a geometrical transformation 
underlying the Hodge theoretic duality away from the stable boundary, i.e. in the bulk of the two moduli spaces involved $\mathcal{M}_{{\rm K3}}$ and $\mathcal{M}_{{\rm het}}^{{\rm G}}$, or at least in the large complex 
structure region \cite{clingher3}. The simplest case to consider is the restriction on the heterotic side to 
the ${\rm P}=0$ locus. From the Hodge-theoretic correspondence, one knows that this 
restriction corresponds on the F-theory side to K3 surfaces with a special lattice polarization of type 
${\rm M}={\rm H} \oplus {\rm E}_8 \oplus {\rm E}_8$. These are precisely the ${\rm M}$-polarized K3 surfaces 
that form the main focus of this paper. Moreover, on the heterotic side, under the vanishing of 
the flat bundle, the B-field has the same properties as a second elliptic curve. Therefore this 
special case of the duality can be regarded as relating Hodge-theoretically 
${\rm M}$-polarized K3 surfaces to pairs of elliptic curves:
$$ {\rm X} \ \longleftrightarrow \ \left ( {\rm E}, \ {\rm B} \right ). $$
This is the precisely the Hodge-theoretic identification from equation $(\ref{hodgecorresp})$. 
\par From this point of view, the transformation we described in Section $\ref{partone}$, and 
for which we have computed explicit formulas in Section $\ref{parttwo}$, provides the proper 
geometrical description of the F-theory/heterotic string duality for ${\rm M}$-polarized K3 surfaces.

\end{document}